\documentclass[twoside,12pt,a4paper]{amsart}
\usepackage{mathrsfs}
\usepackage{amssymb,amsmath,mathrsfs,txfonts,graphicx,color,epsfig}
\usepackage{subfigure}
\usepackage{latexsym}
\usepackage{amsfonts}
\usepackage{amsthm}

\newcommand{\R}{\mathbb{R}}

\newtheorem{theo}{Theorem}[section]
\newtheorem{rem}{Remark}[section]

\newtheorem{lem}{Lemma}[section]

\newtheorem{defi}{Definition}[section]
\numberwithin{equation}{section}

\numberwithin{equation}{section}

\numberwithin{equation}{section}

\begin{document}
\title{\bf Nonlinear structural stability and linear dynamic instability of transonic steady-states
to a hydrodynamic model for semiconductors}

\author{Yue-Hong Feng$^{1,3}$, Ming Mei$^{2,3}$ and Guojing Zhang$^4$}

\date{}

\maketitle \markboth{Y.~H.~Feng,  M.~Mei and G. Zhang} {Structural
stability and dynamical instability of Euler-Poisson equations}

\vspace{-3mm}

\begin{center}
 {\scriptsize
 $^1$College of Mathematics, Faculty of Science, Beijing University of
 Technology, Beijing 100022, China \\[1mm]
$^2$Department of Mathematics, Champlain College Saint-Lambert,
Quebec, J4P 3P2, Canada \\[1mm]
$^3$Department of Mathematics and Statistics, McGill University,
Montreal, Quebec, H3A 2K6, Canada \\[1mm]
$^4$School of Mathematics and Statistics, Northeast Normal University, Changchun 130024,  China \\[1mm]
 Email~: fyh@bjut.edu.cn,\hspace{1mm}
   ming.mei@mcgill.ca, \hspace{1mm} zhanggj100@nenu.edu.cn}
  \end{center}

\

\begin{center}
\begin{minipage}{14cm}
    {\bf Abstract.} {\small
For unipolar hydrodynamic model of semiconductor device represented
by Euler-Poisson equations, when the doping profile is supersonic,
the existence of steady transonic shock solutions and
$C^\infty$-smooth steady transonic solutions for Euler-Poisson Equations
were established
 in \cite{LMZZ18} and
 \cite{WMZZ21},
 respectively. In
this paper  we further study the nonlinear structural stability and
the linear dynamic instability of these steady transonic solutions.
When the $C^1$-smooth transonic steady-states pass through the sonic
line, they produce singularities for the system, and cause some
essential difficulty in the proof of structural stability. For any
relaxation time: $0<\tau\le +\infty$, by means of elaborate
singularity analysis, we first investigate the structural
stability of the $C^1$-smooth transonic steady-states, once the perturbations of the
initial data  and the doping profiles are small enough. Moreover,
when the relaxation time is large enough $\tau\gg 1$, under the condition
that the electric field is positive at the shock location, we prove
that the transonic shock steady-states are  structurally stable with
respect to small perturbations of the supersonic doping profile.
Furthermore, we show the linearly dynamic instability for these
transonic shock steady-states provided that the electric field is
suitable negative. The proofs for the structural stability results
are based on singularity analysis, a monotonicity argument on the
shock position and the downstream density, and the stability
analysis of supersonic and subsonic solutions. The linear dynamic
instability of the steady transonic shock for Euler-Poisson
equations can be transformed to the ill-posedness of a free boundary
problem for the Klein-Gordon  equation. By using a nontrivial
transformation and the shooting method, we prove that the linearized
problem has a transonic shock solution with exponential growths.
These results enrich and develop the existing studies.}
\end{minipage}
\end{center}

\vspace{7mm}


\vspace{5mm}

\noindent {\bf AMS Subject Classification (2000)~:} 35R35, 35Q35,
76N10, 35J70


\vspace{3mm}



\newpage
\tableofcontents

\baselineskip=15pt

\section{Introduction and main results}
\subsection{Modeling equations}

This paper is concerned with the smooth/shock transonic
solutions for the one-dimensional hydrodynamic model for
semiconductors, which is presented as  Euler-Poisson equations with
relaxation effect
\begin{equation}
\label{1.1}
\begin{cases}
    n_t + (nu)_x = 0,   \\
   (nu)_t + \left( P(n) + nu^2 \right)_x = nE - \displaystyle{\frac{nu}{\tau }},  \\
    E_x = n - b(x).
\end{cases}
\end{equation}

This model describes several physical flows including the
propagation of electrons in submicron semiconductor devices
\cite{BDX16,Bl70,DM93,Jungel,MRS90} and plasmas \cite{SM}
(hydrodynamic model), and the biological transport of ions for
channel proteins \cite{CEJS95}. In the hydrodynamic model of
semiconductor devices or plasma, $n$, $u$, $P$ and $E$ represent the
electron density, macroscopic particle velocity, pressure and the
electric field, respectively. The function $b=b(x) > 0$ is the
doping profile standing the impurity for the device.  The parameter
$\tau > 0$ means the momentum relaxation time.
While, the biological model describes the transport of ions between
the extracellular side and the cytoplasmic side of membranes
\cite{CEJS95}. In this case, $n$, $nu$ and $E$ are the ion
concentration, the ions translational mass, and the electric field,
respectively, and the doping profile $b(x)$ represents a background
density of charged ions.

For ideal gas law of isentropic case, the pressure function $P$ is
physically represented by
$$P(n)=Tn^\gamma,$$
where $T  > 0$ is a constant absolute temperature,  and $\gamma > 1$
represents the adiabatic exponent. In this article, we mainly
consider the isothermal case, i.e., $\gamma = 1$ for simplicity of
analysis. The case of $\gamma>1$ can be similarly treated.

Using the terminology from gas dynamics, we call $c := \sqrt{P'(n)}
= \sqrt T
> 0$ the sound speed for $P (n) = T n$. Moreover, if we denote
\begin{equation}
\label{1.3*}
\begin{split}
J := nu,~~\mbox{ the current density of the flow, }  %
 \end{split}
%
\end{equation}
and take $J>0$ without loss of generality, then the flow is
supersonic/sonic/subsonic if the fluid velocity satisfies
\begin{equation}
\label{1.5}
{\rm fluid~~velocity} : u = \frac{J}{n}
> ( {\rm or} ~~  =, ~~  {\rm or} ~~ <)~~c = \sqrt{P'(n)}
= \sqrt T : {\rm sound~~speed}.\end{equation}

Without loss of generality, we assume that $T  = 1$, i.e., $P(n)=
n$. Thus, it can be identified that the flow is subsonic if $n
> J$, sonic if $n =J$, and supersonic if $0 < n < J$.

The main issue of the present paper is to study the structural
stabilities for  $C^1$-smooth transonic steady-states and
transonic shock steady-states, and the linear dynamic instability
for the transonic shock steady-states. These steady-states are
solutions of the following time-independent equations
\begin{equation}
\label{1.2}
\begin{cases}
  (nu)_x = 0, \ i.e., \ J=nu=\mbox{constant},  \\
 \left( n + nu^2 \right)_x = nE - \displaystyle{\frac{nu}{\tau }},  \\
   E_x = n - b(x).
 \end{cases}
\end{equation}
%


For convenience,  we  set
\begin{equation}
\label{1.2*}
\begin{split}
\alpha : = \frac{1}{\tau },~~\mbox{ the  reciprocal of the
relaxation time of
the current},  %
 \end{split}
%
\end{equation}
thus system \eqref{1.2} is reduced to
\begin{equation}\label{1.6}
 \begin{cases}
 \left( n + \displaystyle{\frac{J^2}{n}} \right)_x = nE - \alpha J,   \\
 E_x = n(x) - b(x),
\end{cases}
\end{equation}
equivalently,
\begin{equation}
\begin{cases}\label{new1}
\displaystyle{ n_x=\frac{(nE-\alpha
J)n^2}{n^2-J^2}=\frac{n^2E}{n+J}+\frac{Jn^2(E-\alpha)}{n^2-J^2},
} \\
E_x=n-b.
\end{cases}
\end{equation}

Subjected to the stationary system \eqref{1.6}
 or its equivalent form \eqref{new1},
 two different problems are proposed in this paper.
  One is  the ``initial value" problem for the system \eqref{new1}:
\begin{equation}
\begin{cases}\label{1.7}
\displaystyle{ n_x=\frac{n^2E}{n+J}+\frac{Jn^2(E-\alpha)}{n^2-J^2}},
& x\in \mathbb{R}_+,
 \\
E_x=n-b(x), & x\in \mathbb{R}_+, \\
(n,E)|_{x=0} =(n_0, E_0),
\end{cases}
\end{equation}
in which the initial data $n_0$ is  considered to be   supersonic
satisfying
\[
0<n_0<J
\]
 throughout the paper, because the case of the subsonic data  $n_0>J$ can be similarly treated.

The other is the boundary value problem:
\begin{equation}
\begin{cases}\label{1.4}
\displaystyle{ \left( n + \frac{J^2}{n} \right)_x = nE - \alpha J},
& 0<x<L,
 \\
E_x=n-b(x), & 0<x<L, \\
n(0)=n_l, \ \ n(L) =n_r,
\end{cases}
\end{equation}
where $L>0$, and   the boundary conditions are considered to be
\[
n_l<J<n_r
 \]
 with the supersonic boundary value $n_l$ and  the subsonic boundary value $n_r$ in the paper. The  case of $n_l>J>n_r$ with the subsonic boundary value $n_l$ and  the supersonic boundary value $n_r$ can be also similarly treated.

 The boundary value problem \eqref{1.4}, by dividing \eqref{1.4}$_1$ by $n$ and differentiating it with respect to $x$,  is also equivalent to the following system,
\begin{equation}
\begin{cases}
\displaystyle{ \Big[\Big(\frac{1}{n}-\frac{1}{n^3}\Big)n_x\Big]_x +\alpha \Big(\frac{J}{n}\Big)_x -[n-b(x)]=0,}  & 0<x<L \\
E_x=n-b(x), & 0<x<L, \\
n(0)=n_l, \ \ n(L) =n_r.
\end{cases}
\end{equation}

The $C^1$-smooth transonic solutions and the shock transonic
solutions  for the initial value problem of system \eqref{1.7} or
the boundary value problem \eqref{1.4}  are defined as follows,
respectively.
\begin{defi}[$C^1$-smooth transonic solutions]\label{D1}
A pair of $(n(x), E(x))$ with $n(x)>0$ is called a $C^1$-smooth
transonic solution of the initial value problem \eqref{1.7}, or the
boundary value problem \eqref{1.4}, if $(n,E) \in C^1(I) \times
C^2(I)$ for $I:=\mathbb{R}_+$ or $I:=[0,L]$, and there exists a
number $x_0 >0$ such that
\begin{equation*}
(n, E) =  \left\{ {\begin{array}{*{20}{c}}
  \vspace{3mm} { (n_{\rm sup } , E_{\rm sup } )(x),} \hfill & {{\rm as}\quad x\in [0, x_0],} \hfill  \\
   { (n_{\rm sub}, E_{\rm sub}) (x),} \hfill & {{\rm as}\quad x\ge x_0,} \hfill  \\
\end{array}} \right.
\end{equation*}
where $(n_{sup},E_{sup})(x)$ satisfying $0< n_{\rm sup} (x) < J$ on
$(0, x_0)$ is called to be supersonic, and  $(n_{sub},E_{sub})(x)$
satisfying $n_{\rm sub} (x)
>J$ for $x>x_0$ is subsonic, both of them are differentiable at the sonic line $n=J$:
\begin{equation}\label{9*}
\begin{cases}
n_{\rm sup} (x_0)=n_{\rm sub} (x_0)=J, \quad n'_{\rm sup}
(x_0)=n'_{\rm sub} (x_0),\\
 E_{\rm sup} (x_0)=E_{\rm sub} (x_0), \quad E'_{\rm sup} (x_0)=E'_{\rm sub} (x_0). \quad E''_{\rm sup} (x_0)=E''_{\rm sub} (x_0).
 \end{cases}
\end{equation}
\end{defi}

\begin{defi}[Shock transonic solutions]\label{D2}
A pair of $(n(x), E(x))$ with $n(x)>0$ is called  a transonic shock
solution to problem \eqref{1.7} or \eqref{1.4}, if there exists a
point
 $x_0 \in (0, L)$, such that
\begin{equation*}
\left(  n,E  \right) =  \left\{ {\begin{array}{*{20}{c}}
  \vspace{3mm} {\left( {{n_{\rm sup }},{E_{\rm sup }}} \right)(x),}
  \hfill & {{\rm as}\quad x \in [0, x_0),} \hfill  \\
   {\left( {{n_{\rm sub}},{E_{\rm sub}}} \right)(x),}
    \hfill & {{\rm as}\quad x>x_0,} \hfill  \\
\end{array}} \right.
\end{equation*}
satisfies $0< n_{\rm sup} (x) < J$ on $(0, x_0)$, $n_{\rm sub} (x)
>J$ on $x>x_0$, and  the entropy condition
\begin{equation}\label{Entropy}
0< n_{\rm sup }(x_0^- )< J < n_{\rm sub}(x_0^+ ),
\end{equation}
and the Rankine-Hugoniot conditions
\begin{equation}\label{R-H}
 n_{\rm sup }(x_0^- )
 + \frac{J^2}{n_{\rm sup }(x_0^-
)} = n_{\rm sub}(x_0^+ )  + \frac{J^2}{n_{\rm sub}(x_0^+ )}, \quad
 E_{\rm sup }(x_0^- ) = E_{\rm sub}(x_0^+ ).
\end{equation}
\end{defi}

\subsection{Background of research}

Euler-Poisson equations have been an important topic in fluid
dynamics and semiconductor device industry. One of interesting
questions is to investigate  their physical solutions such as
subsonic/supersonic/transonic solutions. When the setting background
of  the steady-state system of Euler-Poisson equations is completely
subsonic, namely, subsonic boundary and subsonic doping profile,
Degond-Markowich \cite{DM90} first   established the existence of
the subsonic solution,
 and proved its uniqueness once the steady-state system is
 strongly subsonic with $J\ll 1$. Since then, the steady
subsonic flows were studied in great depth with different boundaries
as well as the higher dimensions case in
\cite{BDX14,BDX16,DM90,DM93,FI97,GS05,Je09,Ma91,NS07}, see also the
references therein. These subsonic steady-states with different
settings are then extensively proved to be dynamically stable in
\cite{GS05,LMM,HMW,HMWY1,HMWY2,MWZ,NS07} and the references there
cited, once the initial perturbations around the subsonic
steady-states are small enough.

Regarding the steady supersonic flows, the first result on the
existence and uniqueness of the supersonic steady-states was
obtained by Peng-Violet \cite{PV06}, when the doping profile and the
boundary both are strongly supersonic. See also a recent  study on
supersonic steady-states for 3D potential flows \cite{Bae3}, and the
structural stability of 2-D supersonic steady-states elegantly
proved by Bae-Duan-Xiao-Xie \cite{Bae2}.

As we know, another interesting issue for  Euler-Poisson equations
is about the structure of transonic solutions. The first observation
on such kind transonic shocks was made  by
Ascher-Markowich-Pietra-Schmeiser  \cite{AMPS91} when  the boundary
value problem  \eqref{1.4} is with subsonic boundary data but a
constant supersonic doping profile, which was then generalized by
Rosini  \cite{Ro05} for the non-isentropic flow. Furthermore, when
the doping profile $b(x)$ is nonconstant,  Gamba \cite{Ga92}
constructed 1-D transonic solutions with shocks by  the method of
vanishing viscosity,
 then joint with  Morawetz, they  \cite{GM96} showed the  existence of  transonic solutions with shocks in 2-D case.
 However, these solutions as the limits of vanishing
viscosity yield some boundary layers. Hence, the question of
well-posedness of the boundary value problem for the inviscid
problem couldn't be solved by the vanishing viscosity method. Late
then,  when Euler-Poisson equations  are lack  of the effect of the
semiconductor (the case of $\tau=\infty$), Luo-Xin \cite{LX09}  and
Luo-Rauch-Xie-Xin \cite{LRXX11} studied  the structure of transonic
steady-states, and showed the existence/nonexistence and the
uniqueness/nonuniqueness of the transonic solutions, once the
stationary Euler-Poisson system possess a constant
supersonic/subsonic  doping profile, and one supersonic boundary and
the other
 subsonic boundary. Some
restrictions on the boundary and the domain are also needed. These
transonic shocks with supersonic doping were proved to be
structurally stable \cite{LRXX11}  when the doping profile is a
small perturbation of the constant supersonic doping.  Under certain
restrictions, the time-asymptotic stability of the transonic shock
profiles was also obtained in \cite{LRXX11}.

Recently, the study in this topic has made some profound progress
\cite{CMZZ20,CMZZ21,CMZZ22,LMZZ17,LMZZ18,WMZZ21}.
 For Euler-Poisson equations with relaxation effect
\eqref{1.1},
  when the
boundary is subjected to be sonic (the critical case),
Li-Mei-Zhang-Zhang \cite{LMZZ17,LMZZ18} first classified the
structure of all type of physical solutions. That is, when the
doping profile is subsonic, the steady Euler-Poisson system
possesses a unique subsonic solution, at least one supersonic
solution, and infinitely many shock transonic solutions if the
semiconductor effect is weak ($\tau \gg 1$), and infinitely many
$C^1$-smooth transonic solutions if the semiconductor effect is
strong ($\tau \ll 1$); while, when the doping profile is supersonic
and far from the sonic line, there is no any physical
(subsonic/supersolic/transonic) solution. The supersonic solution
and many shock-transonic solutions exist only when the doping
profile is sufficiently close to the sonic line. Later, when the
doping profile is transonic, according to two cases of the
subsonic-dominated and supersonic-dominated doping profile,
Chen-Mei-Zhang-Zhang \cite{CMZZ20} further classified the structure of all
subsonic/supersonic/shock-transonic solutions.
Very recently, by using the manifold analysis and singularity
analysis near the sonic line and the singular point,
Wei-Mei-Zhang-Zhang \cite{WMZZ21}
 further investigated the existence and regularity of the smooth transonic
steady solutions of Euler-Poisson equations. They gave the detailed
discussions on the structure of directions of the transonic
solutions, and what kind regularity of the smooth transonic
solutions. In particular, when the boundary states are separated in the supersonic regime and
the subsonic regime, they obtained that the Euler-Poisson
system with supersonic doping profile possesses two
$C^\infty$-smooth transonic solutions, where  one is from supersonic region
to subsonic region and the other is of the inverse direction. Moreover, the
existence of  2D and 3D radial subsonic/supersonic/transonic
steady-states with the sonic boundary conditions were
technically proved by Chen-Mei-Zhang-Zhang in \cite{CMZZ21} and
\cite{CMZZ22}, respectively.

%

Remarkably, when the system  \eqref{1.4} is lack of the
semiconductor effect, namely, the relaxation time $\tau=\infty$,
i.e., $\alpha=0$, Luo-Rauch-Xie-Xin \cite{LRXX11} artfully proved
that, there exists a unique transonic shock for the system once the
doping file is a supersonic constant $b=b_0<J$, and showed its
structural stability, namely, there will exist another transonic
shock solution for the  doping profile $b(x)$ as a small
perturbation of the supersonic doping $b_0$. Since the transonic
shocks jump the sonic line from the supersonic regime to the
subsonic regime, such that there is no singularity for the system
\eqref{1.4}  at the sonic line, namely, $n(x)\not=J$. This is an
advantage for the proof of the structural stability of the transonic
shocks, as we know.

However, for the smooth transonic solutions, they pass through the
sonic line $n=J$, and make the system \eqref{new1} to be singular at
the sonic line (the denominator of \eqref{new1} becomes zero, i.e.,
$n^2-J^2=0$). Different from the case of transonic shocks, this
causes an essential difficulty to show the structural stability of
the smooth transonic steady-states, and remains this problem to be
open for any relaxation time $0<\tau\le \infty$. To answer this
question will be one of our main targets in the present paper. Here
we have a key observation. By taking singularity analysis around the
sonic line, we can heuristically determine the value of the
derivative of the smooth transonic solution $n(x)$ at the singular
point on the sonic line. This, with some exquisite singularity
analysis together,  can guarantee us to show the structural
stability around the singular points, then to prove the structural
stability of the smooth transonic steady-state for the initial value
problem  \eqref{1.7} (also the boundary value problem \eqref{1.4})
in the space $C^1_{loc}(\mathbb{R}_+)\times C^2_{loc}(\mathbb{R}_+)$
(or in the space $C^1[0,L]\times C^2[0,L]$). In fact, the
carried-out analysis around the singular transition points on the
sonic line is technical and challenging.

Moreover, when $\tau\gg 1$, we recognize that, the boundary value
problem \eqref{new1} and \eqref{1.4} possess the transonic shocks,
once the doping profile is a supersonic constant, and these
transonic shocks are also structurally stable, when the perturbed
doping profile is small enough. Furthermore, we prove that these
steady transonic shocks are dynamically unstable, when the electric
field is negative. This part can be regarded as the generalizations of
the previous study \cite{LRXX11} with $\tau=\infty$ to the case of
$\tau\gg 1$, but with some technical development.

\subsection{Main results}
In this subsection, we state our main results on the structural
stabilities of smooth transonic steady-states and the transonic
shock steady-states, respectively, and the linear dynamic
instability of these transonic shock steady-states.

We first give the existence and uniqueness of $C^1$-smooth transonic
steady-states and the transonic shock steady-states. This can be
also seen from the following numerical simulations for the phase diagrams of $(n,E)$, for example,
by taking $b(x)=0.5$, $J=1$, and $\tau=1$ in Figures 1 and 2. Here,
there are two smooth curves cross the sonic line $n=J$, namely, two
smooth transonic steady-states.  One smooth transonic curve is from
the supersonic regime to the subsonic regime (see Figure 1) by
setting either the initial data $(n_0,E_0)$ to be supersonic $n_0<J$
or the boundary data to be $n_l<J<n_r$. The other smooth transonic
curve is from the subsonic regime to the supersonic regime (see
Figure 2) by setting either the initial data $(n_0,E_0)$ to be
subsonic $n_0>J$  or the boundary data to be $n_l>J>n_r$.



In what follows, we mainly consider the case of transonic
steady-states from the supersonic regime to the subsonic regime in
Theorems \ref{T1}-\ref{T4}. Of course, the results presented in
Theorem \ref{T1}-\ref{T4} are also true   for the  case of transonic
steady-states from the supersonic regime to the subsonic regime.

\begin{figure}\label{fig0-1}
\centering
\includegraphics[width=0.5\textwidth]{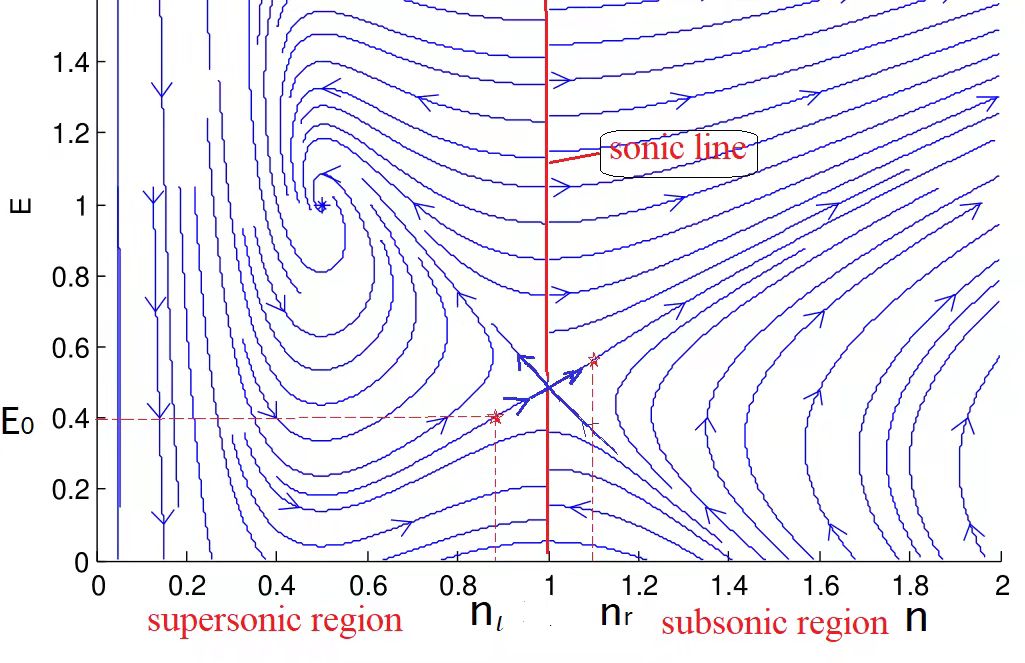}
\caption{This is the $(n,E)$ portrait diagram by taking $b(x)=0.5$, $J=1$, and $\tau=1$. We mark  one smooth curve (smooth transonic steady-state) passing through
the sonic line $n=1$ from the supersonic regime to the subsonic regime.}
\end{figure}

\begin{figure}\label{fig0-2}
\centering
\includegraphics[width=0.5\textwidth]{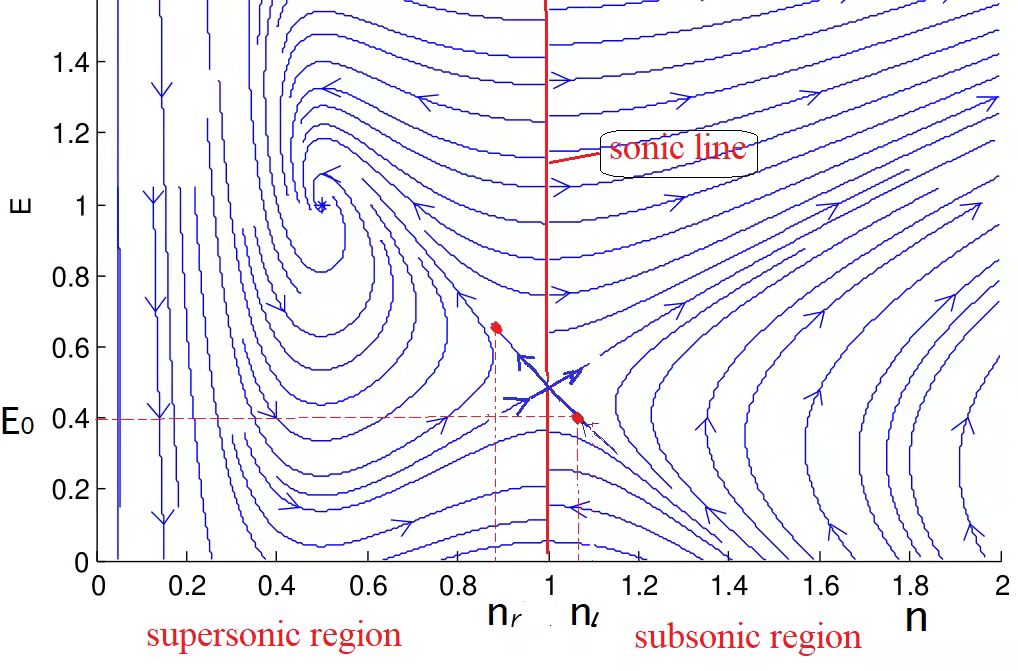}
\caption{This is the $(n,E)$ portrait diagram by taking $b(x)=0.5$, $J=1$, and $\tau=1$. We mark the other smooth curve (smooth transonic steady-state) passing through
the sonic line $n=1$ from the  subsonic regime to the supersonic regime.}
\end{figure}

%
\begin{theo}[Existence and uniqueness  of smooth/shock transonic steady-states]
\label{T1}
Let the doping profile be supersonic such that $b(x) \in L^\infty(0,
L)$ and $0 < b(x) < J$.

\begin{itemize}

\item[(I)] For any relaxation time $0<\tau\le+\infty$, if $b(x)=b_0<J$ is a constant in the supersonic regime, then the
stationary Euler-Poisson equations  with the initial condition
\eqref{1.7}  (or the boundary value condition \eqref{1.4}) admit a
unique $C^1$-smooth transonic solution $(n(x),E(x))$ passing through
the sonic line $n(x)=J$ at a unique point $x_0>0$ determined
implicitly from the system:
\[
(n,E)(x)=\begin{cases}
(n_{sup},E_{sup})(x), & x\in [0,x_0] \\
(n_{sub},E_{sub})(x), & x\ge  x_0,
\end{cases}
\]
and
\[
\begin{cases}
n_{sup}(x_0)=n_{sub}(x_0)=J, \\
  n'_{sup}(x_0)=n'_{sub}(x_0), \\
 E_{sup}(x_0)=E_{sub}(x_0)=\alpha=\frac{1}{\tau}, \\
E'_{sup}(x_0)=E'_{sub}(x_0)=J-b_0, \\
 E''_{sup}(x_0)=E''_{sub}(x_0)=n_x(x_0).
\end{cases}
\]
Here, for the initial value problem \eqref{1.7}, the initial data
$(n_0,E_0)$ is in the supersonic regime, satisfying $0<n_0<J$ and
$E_0<\min\{\frac{1}{\tau},\frac{J}{n_0\tau}\}=\min\{\alpha,\frac{\alpha
J}{n_0}\}$; while for the boundary value problem \eqref{1.4}, the
boundary condition should be suitably selected with $0<n_l<J<n_r$.

\medskip

\item[(II)] If the relaxation time is sufficiently large, $\tau\gg 1$, and
the doping profile $b(x)$ is sufficiently close to the sonic state:
$|b(x)-J|\ll 1$, then the initial value problem \eqref{1.7} (or the
boundary value problem \eqref{1.4}) admits the other type of solution, the so-called transonic shock steady-state $(n,E)(x)$
satisfying the entropy condition
\eqref{Entropy} and the Rankine-Hugoniot jump conditions \eqref{R-H}
at jump location $x_0$,
  which is  unique at  the jump location $x_0$. Here $x_0$ can be uniquely determined when
$n_{\rm sub}(x_0^+ )$ satisfies $\left|n_{\rm sup }(x_0^- )- n_{\rm
sub}(x_0^+ )\right|\ll 1$ is fixed.

\end{itemize}
\end{theo}

\begin{rem}\par

\

\begin{itemize}
\item It is easy to see that the existence results about the $C^1$-smooth
transonic solutions and the transonic shock solutions in Theorem
\ref{T1} are a direct corollary of \cite{WMZZ21} and consequences of
the pioneering works in \cite{LMZZ18}, respectively. So we omit the
details of proof.

\medskip

\item In Part (I), the condition $0<E_0<\min\{\frac{1}{\tau},\frac{J}{n_0\tau}\}=\min\{\alpha,\frac{\alpha J}{n_0}\}$ is to guarantee
$n_{x}(0)>0$ so then the trajectory $(n(x),E(x))$ will pass through
the sonic line with $n=J$ and $E=\alpha$. As showed in
\cite{WMZZ21}, only for certain given initial data $(n_0,E_0)$, the
initial value problem  \eqref{1.7} possesses the unique smooth
transonic steady-state.
\end{itemize}
\end{rem}

Next, we are going to state the structural stability of the
$C^1$-smooth transonic steady-state for the system \eqref{new1} and
\eqref{1.7} as follows.

\begin{theo}[Structural stability of  $C^1$-smooth transonic steady-states of \eqref{1.7}]
\label{T2}
Suppose $J > 0$ to be a constant. For $i=1, 2$, let $b(x)=b_i$ be
two
 constants satisfying $0 <  b_i  < J$ and let $(n_i, E_i)(x)$ be
two $C^1$-smooth transonic solutions (showed in Theorem \ref{T1}) to
the initial value problem \eqref{1.7} with respect to the initial
data $(n_{i0}, E_{i0})$ and  the doping profiles $b_i$, and let  $x
= x_i >0$ be the singular locations of the $C^1$-smooth transonic
solutions $n_i(x)$ cross the sonic line $n(x)=J$, respectively. Then
$(n_i, E_i)(x)$ are structurally stable in
$C^1_{loc}(\mathbb{R}_+)\times C^2_{loc}(\mathbb{R}_+)$. Namely, for
any given local interval $[0,L]\subset \mathbb{R}_+$ with
$L>\max\{x_1,x_2\}$, it holds
\begin{equation}\label{1.10**-2}
 \left\|  \left(n _1 - n _2 \right) \left( \cdot
  \right)
\right\|_{ C^1 [0,L]} +  \left\|  (E_1   -  E_2) (\cdot) \right\|_{
C^2 [0,L] } \le C \delta_0,
\end{equation}
where $C=C(L)>0$ and
\begin{equation}\label{1.10*-2}
\delta_0:=\left| b_1 - b_2 \right| + \left| n _{10} - n _{20}
\right| + \left| E_{10} - E_{20} \right|.
\end{equation}
\end{theo}

Similarly, the structural stability of smooth transonic
steady-states for the boundary value problem \eqref{1.4} holds as
follows.

\begin{theo}[Structural stability of  $C^1$-smooth transonic steady-states of \eqref{1.4}]
\label{T2-2}
Suppose $J > 0$ to be a constant. For $i=1, 2$, let $b_i$ be two
 constants satisfying $0 <  b_i  < J$ and let $(n_i, E_i)(x)$ be
two $C^1$-smooth transonic solutions (showed in Theorem \ref{T1}) to
the boundary problem \eqref{1.4}
 with the boundary data $(n_{il},
n_{ir})$ corresponding to  the doping profiles $b_i$, and let  $x =
x_i >0$ be the singular locations of the $C^1$-smooth transonic
solutions $n_i(x)$ cross the sonic line $n(x)=J$, respectively. Then
$(n_i, E_i)(x)$ for $i=1,2$ are structurally stable in
$C^1[0,L]\times C^2[0,L]$. Namely,  it holds
\begin{equation}\label{1.10**}
 \left\|  \left(n _1 - n _2 \right) \left( \cdot
  \right)
\right\|_{ C^1 [0,L]} +  \left\|  (E_1   -  E_2) (\cdot) \right\|_{
C^2 [0,L] } \le C \delta_0,
\end{equation}
where
\begin{equation}
\delta_0:=\left| b_1 - b_2 \right| + \left| n _{1l} - n _{2l}
\right| + \left| n_{1r} - n_{2r} \right|.
\end{equation}
\end{theo}

\begin{rem} \par

\

\begin{itemize}
\item To our best knowledge, Theorem \ref{T2} and Theorem \ref{T2-2} are the first results to show the structural stability for the smooth transonic steady-states.
\item As showed in \cite{WMZZ21}, the smooth transonic steady-states can be $C^\infty$. By the same fashion as in Theorem \ref{T2} and Theorem \ref{T2-2}, in fact, we can similarly prove their
structural stability in the sense of $C^\infty$. For the simplicity,
we omit it.
\end{itemize}
\end{rem}

Inspired by the study \cite{LRXX11} on the structural stability of steady
transonic shocks for the case with $\tau=\infty$ (i.e., $\alpha=0$)
in \eqref{new1}, we can generalize it to the case with
$\tau\not=\infty$ but $\tau\gg 1$ as follows.

\begin{theo}[Structural stability of transonic shock steady-states of \eqref{1.4}]
\label{T3}
Assume  $J > 0$ is a constant,  the relaxation time is $\tau\gg 1$,
and the doping profile is $0<b(x) = b_0 (x)<J$
 and $|b_0(x)-J|\ll 1$ for $x\in [0, L]$. Let $(n^{(0)},
E^{(0)})(x)$ be the unique transonic shock solution to the boundary
value problem  \eqref{1.4}   with a single transonic shock located
at $x = x_0 \in (0, L)$ satisfying the entropy condition
\eqref{Entropy} and
 the Rankine-Hugoiot condition \eqref{R-H} with $E_{\rm sub}(x_0^+ )>0$.
Then, for a given doping profile $b(x)$ as the small perturbation
around $b_0(x)$, namely, there is $\epsilon_0 > 0$ such that if
\begin{equation}\label{1.10-2}
\left\| b - b_0 \right\|_{C^0[0,L]} = :\epsilon  \le \epsilon _0,
\end{equation}
 the boundary value problem  \eqref{1.4} with $b(x)$ has a
unique transonic shock solution $(\tilde{n}, \tilde{E})(x)$, where
the single transonic shock located at a point $\tilde{x}_0 \in
[x_0-C\epsilon, x_0 + C\epsilon]$ for some constant $C > 0$, namely,
$\tilde{x}_0$ is a small perturbation of $x_0$.
\end{theo}
 The structural stability of the transonic shock steady-states (Theorems \ref{T3}) is also true for the initial value  problem  \eqref{1.7}.

\begin{theo}[Structural stability of transonic shock steady-states of \eqref{1.7}]
\label{T3-2}
Assume  $J > 0$ is a constant,  the relaxation time is $\tau\gg 1$,
and the doping profile is $0<b(x) = b_0 (x)<J$ and $|b_0(x)-J|\ll 1$
for $x\in [0, L]$, where $[0,L]\subset \mathbb{R}_+$ is an any given
subset of $\mathbb{R}_+$. Let $(n^{(0)}, E^{(0)})(x)$ be the unique
transonic shock solution to the initial value problem  \eqref{1.7}
with a single transonic shock located at $x = x_0$ with $0<x_0<L$
satisfying the entropy condition \eqref{Entropy} and
 the Rankine-Hugoiot condition \eqref{R-H} with $E_{\rm sub}(x_0^+ )>0$.
Then, for a given doping profile $b(x)$ as the small perturbation
around $b_0(x)$, namely, there is $\epsilon_0 =\epsilon_0(L)> 0$
such that if
\begin{equation}\label{1.10}
\left\| b - b_0 \right\|_{C^0[0,L]} = :\epsilon  \le \epsilon _0,
\end{equation}
 the initial value problem  \eqref{1.7} with $b(x)$ has a
unique transonic shock solution $(\tilde{n}, \tilde{E})(x)$, where
the single transonic shock located at a point $\tilde{x}_0 \in
[x_0-C\epsilon, x_0 + C\epsilon]$ for some constant $C > 0$, namely,
$\tilde{x}_0$ is a small perturbation of $x_0$.
\end{theo}

Next, we are going to state  the linear dynamic instability of the
steady transonic shock solutions.

For a given function $b(x)$ satisfying $0 < b(x) < J$ for $x \in [0,
L]$, and a constant $\bar J > 0$, let
\begin{equation}\label{1.11}
 ( {\bar n,\bar u,\bar E} )(x) = \left\{
{\begin{array}{*{20}{c}}
   \vspace{3mm}{\left( {{{\bar n}_ - },\displaystyle\frac{{\bar J}}{{{{\bar n}_ - }}},
   {{\bar E}_ - }} \right)(x),} \hfill & {{\rm as} \quad x\in (0, x_0),}
     \\
   {\left( {{{\bar n}_ + },\displaystyle\frac{{\bar J}}{{{{\bar n}_ + }}},
   {{\bar E}_ + }} \right)(x),} \hfill & {{\rm as} \quad x\in (x_0, L),}
\end{array}} \right.
\end{equation}
be a steady transonic shock solution of \eqref{1.2} which satisfies
the boundary conditions
\begin{equation}\label{1.12}
 ( {\bar n,\bar E}  )(0) =  ( {{n_l},{E_l}} ),\quad
\bar n(L) = {n_r},
\end{equation}
where $\bar E(0)=E_{l}$ is determined by the boundary value system
\eqref{1.4}, and  $ ( {\bar n,\bar u,\bar E}  )(x)$ is supersonic as
$ x \in [0, x_0)$, and subsonic as $ x \in (x_0, L]$, i.e.,
\begin{equation}\label{1.13}
\left\{ {\begin{array}{*{20}{c}}
 \vspace{3mm}   \bar n_ - (x)  < \bar J,  \hfill & {{\rm as} \quad
 x\in [0, x_0),} \hfill  \\
   \bar n_ +
 (x) >  \bar J,
   \hfill & {{\rm as} \quad x\in (x_0, L],} \hfill  \\
\end{array}} \right.\end{equation}
and $ (  \bar n,\bar u,\bar E )(x)$ satisfies the Rankine-Hugoniot
conditions at $x = x_0$,
\begin{equation}\label{1.14}
\left( {{\bar n}_ - + \frac{{{{\bar J}^2}}}{{{{\bar n}_ - }}}}
\right)({x_0}) = \left( {{\bar n}_ + + \frac{{{{\bar J}^2}}}{{{{\bar
n}_ + }}}} \right)({x_0}),\quad {\rm and} \quad {\bar E_ - }({x_0})
= {\bar E_ + }({x_0}).\end{equation}
%

Throughout the paper, we also assume that the system is away from
vacuum
\begin{equation}\label{1.15}
\mathop {\inf }\limits_{x \in [0,1]} \bar n(x) > 0.\end{equation}

Obviously, by using the extension Theorem of solutions for  ordinary
differential equations  \cite{Pa92}, we can extend $(\bar n_-, \bar
E_-)$ to be a smooth supersonic solution of \eqref{1.2} on $[0, x_0
+ \delta]$ for some $\delta
> 0$, which coincides with $(\bar n_-, \bar E_-)$ on $[0, x_0]$. In the sequel, we
still use $(\bar n_-, \bar E_-)$ to stand for this extended
solution. In the same way, we shall denote $(\bar n_+, \bar E_+)$ to
be a subsonic solution of \eqref{1.2} on $[x_0-\delta, L]$ for some
$\delta > 0$, which coincides with $(\bar n_+, \bar E_+)$ in
\eqref{1.11} on $[x_0, L]$.

Let us consider the initial boundary value problem of system
\eqref{1.1} with the initial data
\begin{equation}\label{1.16}
\left( {n,u,E} \right)(0,x) = \left( {{n_0},{u_0},{E_0}}
\right)(x),\end{equation}
and the boundary conditions
\begin{equation}\label{1.17}
\left( {n,u,E} \right)(t,0) = \left( {{n_l},\displaystyle
\frac{{\bar J}}{{{n_l}}},{E_l}} \right), \quad n(t,L) =
{n_r},\end{equation}
where $n_l , E_l$ and $n_r$ are the same as that in \eqref{1.12}.

We suppose that the initial values are of the form
\begin{equation}\label{1.18}
\left( {{n_0},{u_0}} \right)(x) = \left\{ {\begin{array}{*{20}{c}}
  \vspace{3mm} {\left( {{n_{0 - }},{u_{0 - }}} \right)(x),}
   \hfill & {{\rm as} \quad x\in (0, {{\tilde x}_0}),} \hfill  \\
   {\left( {{n_{0 + }},{u_{0 + }}} \right)(x),}
   \hfill & {{\rm as} \quad x\in ({{\tilde x}_0},  L),} \hfill  \\
\end{array}} \right.\end{equation}
and
\begin{equation}\label{1.19}
{E_0}(x) = {E_l} + \int_0^x {} \left( {{n_0}(y ) - b(y )} \right)dy
, \end{equation}
which is a small perturbation of $ ( {\bar n,\bar u,\bar E}  )$ in
the sense that
\begin{equation}\label{1.20}
\left|  x_0  -  \tilde x _0  \right| + \left\|  (  n_ {0 +}
 , u_ {0 + }  ) - ( \bar n_ +,\bar u_ +
 ) \right\|_{H^s\left( \left[ \check{x}_0, 1 \right] \right)} +
\left\| ( n_{0 - },u_{0 - } ) -  ( \bar n_ - ,\bar u_ -  )
\right\|_{H^s\left( \left[ 0,\hat x_0 \right] \right)} < \varepsilon
, \end{equation}
for some small $\varepsilon > 0$, and some integer $s$ suitably
large,  where $\check{x}_0 = \min\{x_0, \tilde{x}_0\}$ and
$\hat{x}_0 = \max\{x_0, \tilde{x}_0\}$. Simultaneously, we assume
that $(n_0, u_0, E_0)$  satisfies the Rankine-Hugoniot conditions at
$x = \tilde{x}_0$,
\begin{equation}\label{1.21}
\left( n_{0 + } + n_{0 + }u_{0 + }^2 - n_{0 - } - n_{0 - }u_{0 - }^2
\right)\left( n_{0 + } - n_{0 - } \right)(\tilde x_0) = \left( n_{0
+ }u_{0 + } - n_{0 - }u_{0 - } \right)^2(\tilde x_0).\end{equation}

In advance of declaring our dynamic linear instability results, we
give the definition of piecewise smooth entropy solutions to the
Euler-Poisson equations with relaxation effect \eqref{1.1} as
follows.
\begin{defi}
If $(n_-, u_-, E_-)(x,t)$ and $(n_+, u_+, E_+)(x,t)$ are $C^1$
smooth solutions of Euler-Poisson equations with relaxation effect
\eqref{1.1} in the regions $\{(t, x)|t \ge 0, 0\le x \le s(t)\}$ and
$\{(t, x)|t \ge 0, s(t)\le x \le L\}$, respectively. Then
\begin{equation}\label{1.22}
\left( n,u,E \right)(x,t) = \left\{ {\begin{array}{*{20}{c}}
   \vspace{3mm}{\left( {{n_ - },{u_ - },{E_ - }} \right)(x,t),}
    \hfill & {{  as} \quad x\in (0, s(t)),} \hfill  \\
   {\left( {{n_ + },{u_ + },{E_ + }} \right)(x,t),}
   \hfill & {{  as} \quad x\in (s(t), L),} \hfill  \\
\end{array}} \right.\end{equation}
is said to be a piecewise smooth entropy solution of \eqref{1.1} at
$x = s(t)$ if $(n, u, E)$ satisfies the Rankine-Hugoniot conditions
\begin{equation}\label{1.23}
\begin{cases}
    \left(  n +  n u^2  \right)(t,s(t)^+ )
 - \left( n + n u^2 \right)(t,s(t)^- ) = \left( nu(t,s(t)^+ )
  - nu(t,s(t)^- ) \right)  s'(t),    \\
     \left( nu(t,s(t)^+ ) - nu(t,s(t)^- )  \right)
  = \left(  n(t,s(t)^+ ) - n(t,s(t)^- )  \right)  s'(t),  \\
    E(t,s(t)^+ ) = E(t,s(t)^- ),
\end{cases}
\end{equation}
and the Lax geometric entropy condition
\begin{equation*}
\left(  u  - 1   \right)(t,s(t)^- )  >  s'(t)  >  \left( u   -   1
\right)(t,s(t)^+ ), ~~{\rm and} ~~ \left(  u   +  1 \right)(t,s(t)^+
)  >  s'(t).
\end{equation*}
\end{defi}

Now the linear dynamic instability theorem in this paper is declared
as follows.

%
\begin{theo}[Linearly dynamic instability of transonic shock steady-states]
\label{T4}
Let $(\bar n, \bar u, \bar E)(x)$ be a  transonic shock steady-state
to system \eqref{1.1} satisfying \eqref{1.11}-\eqref{1.15}. There
exists
 $\delta > 0$ such that if
\begin{equation}\label{1.24}
{\bar E_ - }({x_0}) = {\bar E_ + }({x_0}) <  - \delta,
\end{equation}
then the linearized problem corresponding to the initial boundary
problem \eqref{1.1} and \eqref{1.16}-\eqref{1.21} admits a linearly
unstable transonic shock solution $( n, u,  E)(x,t)$ which is
time-exponentially growing away from the  transonic shock
steady-state $(\bar n, \bar u, \bar E)(x)$.
\end{theo}
\begin{rem}
There is a fundamental difficulty that the problem involves a free
boundary (shock) on the left of the subsonic region. To overcome
this embarrassment, the key idea is to introduce a nontrivial
transformation to reformulate the problem on the fixed domain $[x_0,
L]$.\end{rem}

\subsection{Strategies for proofs}

In this subsection, we are going to state the ideas and strategies
for proving Theorems \ref{T2}-\ref{T4}.

For the structural stability of the smooth transonic steady-states
 stated in Theorem \ref{T2} and Theorem \ref{T2-2},
the key steps are to carry out the singular analysis around
 the singular points when the smooth transonic steady-states
 cross the sonic line.
Since there are some singularities for Euler-Poisson
  equation \eqref{new1} around the singular point,
  a suitable setting and re-organizing for the working
  system in $[J-\varepsilon_*,J+\varepsilon_*]$ are quite
  technical and artful. The total procedure of proof
   will be divided in  two cases $\alpha=0$ (i.e., $\tau=\infty$) and
 $\alpha>0$ (i.e., $\tau\not=\infty$),  and use six steps represented by six
   lemmas (see Lemmas \ref{L2-1}-\ref{L2-6}).
We first treat the easy case of $\alpha=0$. In fact, when $\alpha=0$, problem \eqref{1.7} is reduced to a variable
separable  ordinary differential equation, and then we get the
explicit formula of the corresponding trajectory (see \eqref{1-8}
for $E (n) = (n - J)W(n,b)$). Furthermore, with the help of the
property on $W(n,b)$ obtained in Lemma \ref{L2-1}, we successfully
overcome the difficulty caused by the singularity and obtain the
stability result in the first case of Theorem  \ref{T2}.
When $\alpha>0$, the task becomes very difficult. Different from the
case of $\alpha=0$, there is no any explicit formula for $\tilde
W(n,b)$ since problem \eqref{1.7} can't turn into a variable
separable ordinary differential equation.
We first introduce a transformation $ \tilde E  =E - \frac{ \alpha J
}{n}$, and then get the corresponding trajectory equation (see
\eqref{2-4}) to the reduced problem for  \eqref{1.7}. After that, we
study the properties of $\frac{\tilde E}{n-J}=\tilde W(n,b)$ with
respect to variable $n$ and parameter $b$ in Lemma \ref{L2-4} and
Lemma \ref{L2-5}, respectively. Next, we set $M=n(L)$ and translate
the domain $x\in [0, L]$ into the targeted domain $n\in [n_*, M]$.
Late then, we split the targeted  domain $[n_*,M]$ into three
 parts $[n_*, J-\varepsilon_*]\cup[J-\varepsilon_*,J+\varepsilon_*]
\cup [J+\varepsilon_*, M]$, where
$[J-\varepsilon_*,J+\varepsilon_*]$ is the singular domain including
the
  singular points $n=J$,
  and $[n_*, J-\varepsilon_*]
\cup [J+\varepsilon_*, M]$ are the non-singular domains. The crucial
process is to evaluate  the difference of two smooth transonic
steady-states $|(n_1-n_2)(x)|+|(E_1-E_2)(x)|$
 and $|(n_1-n_2)_x(x)|+|(E_1-E_2)_{xx}(x)|$ near
  the singular point in $[J-\varepsilon_*,J+\varepsilon_*]$.
We use the difference scheme and the manifold analysis near the
singularity point $n=J$ to remove the singular property of $\tilde
W(n,b)$. By the method of proof by contradiction, we can fix a
positive constant $\varepsilon_*=\min\{\varepsilon_+,
\varepsilon_-\}$ suitably small, and prove that $\frac{\mathcal
{P}}{n-J}$ admits both the upper bound
 $l-\delta_*$ and the lower bound  $l+\delta_*$  on the domain
 $n \in
[J-\varepsilon_*,J+\varepsilon_*].$ Next, due to the fact that there
is no singularity on the domain $[n_*, J-\varepsilon_*]\cup
[J+\varepsilon_*, M]$, we easily obtain $\left|\frac{\mathcal
{P}}{n-J}\right|<C$ over  the targeted  domain $ n\in [n_*, M].$
Furthermore, by combining the well-established estimates, we obtain that
$\tilde W(n,b)$ is Lipschitz continuous with respect to the
parameter $b$. Finally, by combining Lemmas \ref{L2-3}-\ref{L2-5},
we prove Lemma \ref{L2-6} which contains the structural stability
 of  $C^1$-smooth transonic steady-states of \eqref{1.7} for the second
 case of Theorem \ref{T2}.

For the structural stability of steady transonic shocks stated in Theorem
\ref{T3} and Theorem \ref{T3-2}, the main idea is based on  a
monotonic dependence of the shock location as a function of
downstream density and a priori estimates for supersonic and
subsonic solutions.
First, by the entropy condition and the Rankine-Hugoniot condition,
we connect a supersonic state $(n,E)$ satisfying $n < J$ to a unique
subsonic state $(\mathscr{S}(n);E)$ via a transonic shock. Next,
with the help of the positive electric field condition $E_{\rm sup
}^{(2)}({x_1}) > 0$ and the comparison principles for  ordinary
differential equations, we establish the monotonic relation for the
transonic shock solutions (see Lemma \ref{L2.1}). Then, by using the
multiplier method, we establish the a priori estimates for
supersonic and subsonic flows, which yield the existence of
supersonic, subsonic, and transonic shock solutions (see Lemma
\ref{L2.2}). After that, we start to prove Theorem \ref{T3}. Based
on the fact that the boundary value problem  \eqref{1.4} has a
unique transonic shock solution $(n^{(0)}, E^{(0)})$ for the case
when $b(x) = b_0(x) (x \in [0, L])$ with a single transonic shock
located at $x = x_0 \in (0, L)$, we construct two different
transonic shock solutions whose subsonic solutions $(n^r_i,
E^r_i)(x)$, $(i=1,2)$ on the interval $x \in [x_i , L]$, in which
shock locations are $x_1=x_0-\delta$ and $x_2=x_0+\delta$,
respectively. Therefore, it follows from Lemma \ref{L2.1} that
$n_2^r(L) < {n_r} < n_1^r(L)$. Late then, based on $\left(
{n_1^r,E_1^r} \right)$ and $\left( {n_2^r,E_2^r} \right)$, we
further define two transonic solutions $( {{{\hat n}^{(i)}},{{\hat
E}^{(i)}}} )(x)$, $(i=1,2)$, as $b$ is a small perturbation of
$b_0$. And then, by Lemma \ref{L2.2}, we obtain $\left|\hat n_i^r(L)
- n_i^r( L )\right| \le C\varepsilon$, $(i=1,2)$. Finally, the
desiring stability result in Theorem \ref{T3} follows by combining
the above estimates and a monotonicity argument. We find that the
boundary problem \eqref{1.4} admits a unique transonic shock
solution $(\tilde n, \tilde E  )$ with a single transonic shock
located at some point $\tilde{x}_0 \in (x_1, x_2)$.

Now, let us explain  the key difference between the proofs of
Theorem \ref{T2} (similarly Theorem \ref{T2-2}) and  Theorem
\ref{T3} (similarly Theorem \ref{T3-2}).  Since the solutions
considered in Theorem \ref{T3} are the transonic shocks, which jump
from the supersonic region to the subsonic region, and do not
directly cross the sonic line. So there is no singularity for the
system near the sonic line. This is a kind of advantage in the proof
of structural stability. However, for Theorem \ref{T2} and Theorem
\ref{T2-2},  the smooth transonic steady-states pass through the
sonic line, which cause the working system \eqref{new1} to be
singular. This is essentially different and also challenging in the
proofs.

In what follows, we talk about the strategy for the proof of the
linear dynamic instability in Theorem \ref{T4}. Although the idea
comes from the previous study in \cite{LRXX11}, it is still not
straightforward. First, by the Rankine-Hugoniot conditions and the
implicit function Theorem, we formulate an initial boundary value
problem in the region $\{(t, x)| t>0, x> s(t)\}$. Next, we introduce
a nontrivial transformation to reformulate this free boundary
problem into a fixed  boundary problem.  After that, we get the
linearized initial boundary value problem \eqref{3.23} for
consideration. Hence, in view of problem \eqref{3.23} resembles a
Klein-Gordon equation, we prove that it admits a transonic shock
solution with exponential
growths by  the shooting method.

We end this section by stating the arrangement of the rest of this
paper. In Section 2, we establish the structural stability for the
steady $C^1$-smooth transonic solution, by carrying out the singular
analysis near the sonic line.
 In Section 3, we
show the structural stability of the steady transonic shock
solutions. We
 first give two useful lemmas which include the monotonic
relations for the transonic shock solutions and the a priori
estimates for supersonic and subsonic flows. Then, we use three
steps to complete the proof of the Theorem \ref{T3}.
 In the last section, we
study the linear dynamic instability of transonic shock solutions.
We  formulate the linearized problem, and then construct a shock
solution with exponential growths to complete the proof of Theorem
\ref{T4}.

%
%
\section{Structural stability for steady $C^1$-smooth transonic solutions}\label{s2}

This section is to devoted to the proof of structural stability of
smooth transonic steady-states stated in Theorem \ref{T2} and
Theorem \ref{T2-2}. Here we mainly give the detailed proof
 to Theorem \ref{T2}, because Theorem \ref{T2-2} can be similarly done. The proof is divided into two cases: $\alpha=\frac{1}{\tau}=0$ and $\alpha=\frac{1}{\tau}>0$.
 We first investigate the structural stability of the smooth transonic steady-states in the easy case of $\alpha=0$. The advantage in this case is that the electric field $E=E(n)$ can be explicitly expressed, which makes the singularity analysis to be simple and direct, and can help us to build up the structural stability. Secondly, we treat the case of $\alpha>0$. Since the relationship of $E=E(n)$ is implicit, the singularity for the system of equations for $(n,E)$ crossing the sonic line $n=J$ causes us an essential difficulty. So, some technical analysis around the singular points needs to be artfully carried out. This will be the crucial step for the proof of the structural stability of the smooth transonic steady-states.

\subsection{Case 1.  $\alpha=0$ (i.e.,  $\tau=+\infty$)}
For $\alpha=0$, problem \eqref{1.7} becomes
\begin{equation}\label{1-1}
 \begin{cases}
 \left( n + \displaystyle{\frac{J^2}{n}} \right)_x = nE,   \\
 E_x = n(x) - b, \\
(n,E)|_{x=0} =(n_0, E_0).
\end{cases}
\end{equation}
When $n\neq J$, problem \eqref{1-1} is equivalent to
\begin{equation}
\begin{cases}\label{1-2}
\displaystyle{ n_x=\frac{n^3E}{n^2-J^2},}\\
E_x=n-b, \\
(n,E)|_{x=0} =(n_0, E_0).
\end{cases}
\end{equation}
Then the trajectory for the equations of problem \eqref{1-2} is
$$EdE =
\frac{\left( n^2 - J^2 \right)(n - b)}{n^3}dn.$$
Integrating it, we have
$$\frac{1}{2}{E^2} = n - b\ln n + \frac{ J^2 }{n} -
\frac{ b J^2 }{ 2 n^2 } + C_0,$$
where $C_0$ is a constant to be determined by
\[
C_0=b\ln J +\frac{b}{2}-2J,
\]
due to the fact that the curve of the $C^1$-smooth transonic
solution must pass through the point $(n,E)=(J,0)$.

\begin{figure}
\centering
\includegraphics[width=6cm]{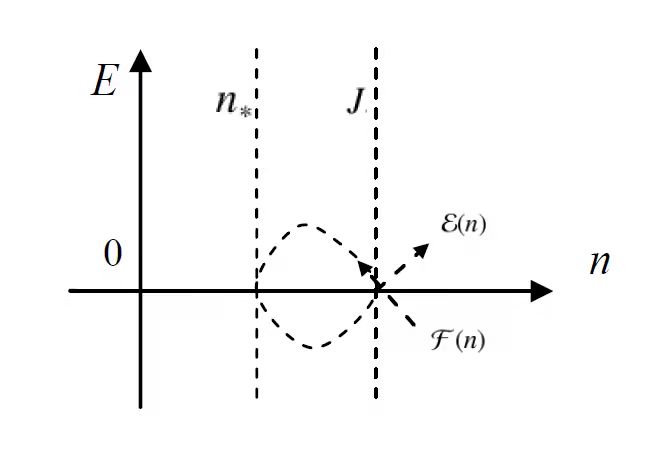}
\caption{The trajectory of $E$ and $n$.}
\end{figure}


Thus, we have
\begin{equation}
\label{1-3}
E^2 = 2n - 2b\ln n + \frac{2J^2}{n} - \frac{bJ^2}{n^2} + 2C_0
\buildrel \Delta \over = g(n),
\end{equation}
with $g(n_*)=0$, where $0<n_*<n<+\infty$.

Let us use $\mathcal {E}(n)$ to denote the trajectory which is from
the supersonic region to the subsonic region, and use $\mathcal
{F}(n)$ to stand for the other which has a reverse direction (see
Figure 3). In the following, we only consider the $C^1$-smooth
transonic solution corresponding to $\mathcal {E}(n)$.
For this end, we want to investigate the properties of $\mathcal
{E}(n)$ to remove the singularity of the targeted equations.

Obviously, from \eqref{1-3}, we have
\begin{equation}
\label{1-4}
 E (n) = \left\{  \begin{split}
   - \sqrt {g(n)},   \quad  n \le J,   \\
   \sqrt {g(n)  },  \quad  n \ge J.
\end{split} \right.
\end{equation}
Here $g(n)$ satisfies $g(J)=0$ and
\begin{equation}
\label{1-5}
\left. \frac{dg}{dn} \right|_{n = J} = \left( 2 - \frac{2b}{n} -
\frac{2J^2}{n^2} + 2\frac{bJ^2}{n^3} \right)_{n = J} = 0,
\end{equation}
%
%
\begin{equation}
\label{1-6}
\left. \frac{d^2g}{dn^2} \right|_{n = J} = \left( \frac{2b}{n^2} +
\frac{4J^2}{n^3} - 6\frac{bJ^2}{n^4} \right)_{n = J} =
\frac{4}{J}\left( 1 - \frac{b}{J} \right) > 0,
\end{equation}
and
\begin{equation}
\label{1-7}
\frac{d^3g}{dn^3}(n) =  - \frac{4}{n^3}\left( b + \frac{3J^2}{n} -
\frac{6bJ^2}{n^2} \right) \buildrel \Delta \over = h(n,b),
\end{equation}
where we have used the supersonic doping profile condition $b<J$.

Then it follows from the facts $g(n)\in C^\infty (n_*, \infty)$,  \eqref{1-5}-\eqref{1-7}
and the Taylor's formula with the integral remainder,  that
\begin{equation*}
\begin{split}
 g(n) =\ & g(J) + (n - J) \left.  \frac{ dg }{dn} \right|_{n = J}
  + \frac{(n - J)^2}{2}\left. \frac{d^2g}{dn^2} \right|_{n = J}
   + \frac{(n - J)^3}{2}\int_0^1  (1 - t)^2
   \frac{d^3g}{dn^3}(J + t(n - J))dt \\
  = \ &\frac{(n - J)^2}{2}\frac{4}{J}\left( 1 - \frac{b}{J} \right)
   + \frac{(n - J)^3}{2}\int_0^1  (1 - t)^2
   \frac{d^3g}{dn^3}(J + t(n - J))dt \\
  = \ &(n - J)^2\left( \frac{2}{J}\left( 1 - \frac{b}{J} \right)
  + \frac{n - J}{2}\int_0^1 (1 - t)^2h(J + t(n - J),b)dt \right),
 \end{split}
\end{equation*}
which implies
\begin{equation}
\label{1-8}
 E (n)  \buildrel \Delta \over = (n - J)W(n,b),
\end{equation}
where
\begin{equation}
W(n,b)= \sqrt {\frac{2}{J}\left( 1 - \frac{b}{J} \right) +
\frac{n - J}{2}\int_0^1  (1 - t)^2h(J + t(n - J),b)dt}.
\end{equation}

It is easy to see that $W(n,b)$ has the following properties.

\begin{lem}\label{L2-1}
For $M>n_*$ and $0<\varepsilon<J$, there exist positive constants
$C_1$ and $C_2$ depend only on $\varepsilon$ and $M$ such that
\begin{equation}\label{1-9}
  0<C_1<W(n,b)<C_2, \quad for ~~n\in(n_*, M)~~ and~~ b\in(\varepsilon,
  J).
    \end{equation}
Moreover, it holds
\begin{equation*}
 W(n,b)\in C^\infty((n_*,M)\times (\varepsilon, J)),
    \end{equation*}
and
\begin{equation*}
 |\partial_n W(n,b)|<C, \quad |\partial_b W(n,b)|<C, \quad for
 ~~n\in (n_*,M),
 ~~ and ~~b\in (\varepsilon, J).
    \end{equation*}
\end{lem}

Now, we are going to investigate the structural stability of the $C^1$-smooth
transonic steady-states.

For $i=1,2$, let $(n_i, E_i)$ be the two $C^1$-smooth transonic steady-states
satisfying
\begin{equation}\label{1-10}
  \left\{  \begin{split}
&  n_{ix}  = \frac{n_i^3}{n_i + J}\frac{E_i}{n_i - J}
=\frac{n_i^3}{n_i + J}W(n_i, b_i),  \\
& E_{ix} = n_i - b_i,\quad x \in [0,L], \\
& n_i(0) = n_{i0},\quad E_i(0) = E_{i0},
\end{split} \right.
    \end{equation}
where $0<n_{i0}<J, E_{i0}=\mathcal {E}(n_{i0})$.

\begin{lem}\label{L2-2}
There exists a constant $C>0$  such that
\begin{equation}\label{1-11}
\|  n _1(\cdot) -  n _2 (\cdot)\|_{ C^1 [0,L]}
  + \|  E_1(\cdot) -  E_2 (\cdot)\|_{ C^2 [0,L]} \le C\delta_0,
    \end{equation}
where $\delta_0$ is the same meaning as that in \eqref{1.10*-2}.
\end{lem}
\noindent\textbf{Proof.}
By making difference of \eqref{1-10} with respect to $n_1$ and $n_2$, we get
\begin{equation*}\left\{ \begin{split}
  &  \left( n_1 - n_2 \right)_x = \frac{n_1^3}{n_1 + J}W(n_1,b_1)
   - \frac{n_2^3}{n_2 + J}W(n_2,b_2),   \\
  & \left( n_1 - n_2 \right)(0) = n_{10} - n_{20},\quad x \in [0,L].
\end{split} \right.\end{equation*}
Then, by Lemma \ref{L2-1}, we have
\begin{equation}\label{1-12}
 \begin{split}
 \left( n_1 - n_2 \right)_x =&
  \left( \frac{n_1^3}{n_1 + J}
   - \frac{n_2^3}{n_2 + J} \right)
   W(n_1,b_1) + \frac{n_2^3}{n_2 + J}\left( W(n_1,b_1) - W(n_2,b_2) \right) \\
  \le & C\left| n_1 - n_2 \right| + C\left| b_1 - b_2 \right|,
 \end{split}
    \end{equation}
which implies
\begin{equation*}
\frac{d\left( n_1 - n_2 \right)^2}{dx} \le C\left( \left| n_1 - n_2
\right|^2 + \left|b_1 - b_2 \right|^2 \right), 
\end{equation*}
where the Cauchy-Schwarz inequality was used. Following the same
way, we find that \eqref{1-12} is also true for $n_2-n_1$.
%
%
Therefore, it follows that
\begin{equation}\label{1-13}
 \begin{split}
\left( n_1 - n_2 \right)^2 \le C\left( \left| n_{10} - n_{20}
\right|^2 + \left| b_1 - b_2 \right|^2 \right), \quad x\in [0,L],
 \end{split}
    \end{equation}
namely
\begin{equation}
\label{1-14}
 \left| n_1 - n_2 \right| \le C\left( \left| n_{10} - n_{20} \right|
+ \left| b_1 - b_2 \right|  \right)\quad x\in [0,L].
    \end{equation}
This together with \eqref{1-12} yields
\begin{equation} \label{1-15}
\left| (n_1 - n_2)_x \right| \le C\left( \left| n_{10} - n_{20}
\right| + \left| b_1 - b_2 \right|  \right),\quad x\in [0,L].
    \end{equation}
On the other hand, from the second equation of \eqref{1-10}, we
obtain
\begin{equation} \label{1-16}
  E_i (x) =  E_{i0}  + \int_0^x
 \left(  n _i (y) - b_i  \right)dy, \quad i=1,2. \end{equation}
Then it follows  that
\begin{align}\label{1-17}
   \left| {{E_1}(x) - {E_2}(x)} \right| \le& \left| {{E_{10}} - {E_{20}}} \right| + \int_0^x {} \left| {{n _1}(y) - {n _2}(y)} \right|dy + x\left| {{b_1} - {b_2}} \right| \\
  \le &\left| {{E_{10}} - {E_{20}}} \right|
  + C\mathop {\sup }\limits_{y \in [0,L]}
  \left| {{n _1}(y) - {n _2}(y)} \right| + \left| {{b_1} - {b_2}} \right|\notag \\
  \le &C\delta_0,\quad x\in [0,L],\notag \end{align}
\begin{equation}\label{1-18}
 |E_{1x}(x)-E_{2x}(x)|=|n_{1}(x)-n_{2}(x)-(b_1-b_2)|\le C \delta_0,
,\quad x\in [0,L],
\end{equation}
and
\begin{equation}\label{1-19}
|E_{1xx}(x)-E_{2xx}(x)|=|n_{1x}(x)-n_{2x}(x)|\le C \delta_0,\quad
x\in [0,L].
\end{equation}
Therefore, by combining \eqref{1-14}- \eqref{1-15} and \eqref{1-17}-
\eqref{1-19}, we obtain
\[
  \|  n _1(\cdot) -  n _2 (\cdot)\|_{ C^1 [0,L]}
  + \|  E_1(\cdot) -  E_2 (\cdot)\|_{ C^2 [0,L]} \le C\delta_0.
 \]
The proof of Lemma \ref{L2-2} is completed. \hfill $\Box$

\subsection{Case 2. $\alpha> 0$ (i.e.,  $0<\tau<+\infty$)}
In this subsection, we continue to study the following problem with $\alpha>0$
\begin{equation}
\begin{cases}\label{2-0}
 \left( n + \displaystyle{\frac{J^2}{n}} \right)_x = nE - \alpha J,
& x\in [0,L],
 \\
E_x=n-b, & x\in [0,L], \\
(n,E)|_{x=0} =(n_0, E_0).
\end{cases}
\end{equation}
Obviously, the first equation of \eqref{2-0} can be rewritten as
\begin{equation}\label{2-1}
\frac{(n - J)(n + J)}{n^2}n_x = nE - \alpha J \buildrel \Delta \over
= n\tilde E,
\end{equation}
where
$$ \tilde E =\tilde E (n,b) :=E - \frac{ \alpha J }{n}.
$$
In view of \eqref{2-1} and %
\begin{equation}\label{2-2}
\tilde E_x = E_x + \frac{\alpha J}{n^2}n_x = n - b + \frac{\alpha
J}{n^2}n_x,
\end{equation}
it follows that the unknowns $(n, \tilde E)$ satisfies
\begin{equation}\label{2-3}
\left\{ \begin{split}
& n_x = \frac{n^3}{n + J}\frac{\tilde E}{n - J}, \\
& \tilde E_x = n - b + \frac{\alpha J}{n^2}n_x.
 \end{split} \right.
\end{equation}
Then  the corresponding trajectory equation to system \eqref{2-3} is
\begin{equation}\label{2-4}
  \frac{d\tilde E}{dn} = \frac{(n + J)(n -
b)(n - J)}{n^3\tilde E} + \frac{\alpha J}{n^2} = \frac{(n + J)(n -
b)}{n^3}\frac{1}{\frac{\tilde E}{n - J}} + \frac{\alpha J}{n^2}.
\end{equation}

It follows from \cite{WMZZ21} that  Euler-Poisson system
\eqref{2-3} or \eqref{2-4} possesses two $C^\infty$-smooth transonic
solutions. One is denoted by $\tilde{\mathcal {E}}(n)$ which is from
supersonic region to subsonic region, and the other is
$\tilde{\mathcal {F}}(n)$ which has the inverse direction.
   In the
following, we only consider the $C^1$-smooth transonic solution
 $\tilde{\mathcal {E}}(n)$.
Let $\tilde E= \tilde E(n)= \tilde E(n,b)$ be the trajectory
corresponding to $\tilde{\mathcal {E}}(n)$. Then from \cite{WMZZ21},
the property of $\tilde E(n,b)$ is stated as follows.
\begin{lem}\label{L2-3}
 $\tilde E(n,b)$  is smooth respect to $n$, and is continuous about
 $b$. Namely,
  $$\tilde E(n,b)\in C^\infty(n_*, \infty)\times C^0(0,J),$$
where $n_*$ satisfies $\tilde E(n_*,b)=0$ and  $0<n_*<J$.
\end{lem}

Set%
 $$\lim  \limits_{n \to J} \frac{ d\tilde E(n) }{ dn } = k.$$
By \eqref{2-4} and the Hospital's rule, we have
$$k = \frac{2(J - b)}{J^2}\frac{1}{k} + \frac{\alpha }{J},$$
which implies that
\begin{equation}\label{2-5}
 k = k_ \pm  = \frac{1}{2}\left( \frac{\alpha }{J} \pm \sqrt {\left(
\frac{\alpha }{J} \right)^2 + \frac{8(J - b)}{J^2} } \right).
\end{equation}
Since the targeted trajectory is $\tilde{\mathcal {E}}(n)$ and
$b<J$, we get $k = k_ +>0.$

In order to prove the stability of $C^1$-smooth transonic solution,
the analysis of properties of $\frac{\tilde E}{n-J}$ is crucial due
to \eqref{2-3}. Particularly, we have to investigate the property of
$\frac{\tilde E}{n-J}$ about the parameter $b$.
In what follows, let us denote $$\frac{\tilde E}{n-J}=\tilde
W(n,b).$$ Then, in the similar fashion in \cite{WMZZ21}, we have

\begin{lem}\label{L2-4}
For $0<b<J$, $\tilde W(n,b)$  is smooth respect to $n$, and is
continuous about
 $b$. Moreover, there exist constants $M=n(L)>n_*$ and $C=C(M,b)>0$ such that
  $$|\tilde W(n,b)|<C,\quad \forall n\in  [n_*, M], $$
and $$\left|\frac{d\tilde W(n,b)}{dn}\right|<C,\quad \forall n\in
[n_*, M],$$
 where $n_*$ is the same meaning as that in Lemma
\ref{L2-3}.
\end{lem}
\begin{rem}
To prove the property of $\tilde W(n,b)$ about the parameter $b$ is
key but difficult. In the case of $\alpha=\frac{1}{\tau}=0$, the
proof of the property about $W(n,b)$ is easy since it has an
explicit formula. However,  $\tilde W(n,b)$ doesn't has the explicit
representation due to the fact that \eqref{2-4} is not separated
type.
\end{rem}

Next, we begin to establish some necessary estimates and to prove that
 $\tilde W(n,b)$  is Lipschitz continuous with respect to the parameter $b$ as follows,
\begin{lem}\label{L2-5}
For $0<b_1, b_2<J$, there exists constant $C>0$ such that
\begin{equation}\label{2-6}
 \left| \tilde W(n,b_1)- \tilde W(n,b_2)\right|\le C|b_1-b_2|, \quad \forall
 n\in [n_*, M].
  \end{equation}
\end{lem}
\noindent\textbf{Proof.} For $h>0$, by choosing $b_1=b+h$ and
$b_2=b$, it follows from \eqref{2-4} that
$$\frac{d\tilde E(n,b + h)}{dn} = \frac{(n + J)\left( n - (b + h)
\right)(n - J)}{n^3\tilde E(n,b + h)} + \frac{\alpha J}{n^2},$$
and
$$\frac{d\tilde E(n,b)}{dn} = \frac{(n + J)\left( n - b
\right)(n - J)}{n^3\tilde E(n,b )} + \frac{\alpha J}{n^2}.$$

By taking difference of the above two equations, we have
\begin{eqnarray}\label{2-7}
&&\frac{d}{dn} ( {\tilde E}(n,b + h) - {\tilde E}(n,b) )  \\
&&  = \left( \frac{(n + J)\left( n - (b + h) \right)
  (n - J)}{n^3\tilde E(n,b + h)} - \frac{(n + J)
  \left( n - b \right)(n - J)}{n^3\tilde E(n,b + h)} \right) \notag \\
& & \ \ \ + \left( \frac{(n + J)\left( n - b \right)(n - J)}{n^3
  \tilde E(n,b + h)} - \frac{(n + J)
  \left( n - b \right)(n - J)}{n^3\tilde E(n,b)} \right) \notag \\
&&  = \frac{n + J}{n^3}\frac{1}{\frac{\tilde E(n,b + h)}{n - J}}
  \left(  - h \right) + \frac{(n + J)\left( n - b \right)}{n^3}
  \left( \frac{n - J}{\tilde E(n,b + h)} - \frac{n - J}{\tilde E(n,b)}
  \right)  \notag\\
 && = \frac{n + J}{n^3}\frac{1}{\frac{\tilde E(n,b + h)}{n - J}}
  \left(  - h \right) +  \frac{(n + J)\left( n - b \right)}{n^3}
  \frac{(n - J)^2}{\tilde E(n,b + h)\tilde E(n,b)}
  \frac{\tilde E(n,b) - \tilde E(n,b + h)}{n - J} \notag \\
&&  = \frac{n + J}{n^3}\frac{1}{\frac{\tilde E(n,b + h)}{n - J}}
  \left(  - h \right) +  \frac{(n + J)\left( n - b \right)}{n^3}
  \frac{1}{\frac{\tilde E(n,b + h)}{n - J}\frac{\tilde E(n,b)}{n - J}}
  \frac{\tilde E(n,b) - \tilde E(n,b + h)}{n - J}. \notag
  \end{eqnarray}
 Dividing \eqref{2-7} by $h$, and letting
 $$\frac{\tilde E(n,b + h) - \tilde E(n,b)}{h}
 \buildrel \Delta \over = \Delta _h^b \tilde E(n,b),
$$
we obtain
$$\frac{d\left( \Delta _h^b\tilde E(n,b) \right)}{dn} \buildrel
\Delta \over = f(n,b,h) - g(n,b,h)\frac{\Delta _h^b\tilde E(n,b)}{n
- J},$$
where
$$f(n,b,h) =  - \frac{n + J}{n^3}\frac{1}{\frac{\tilde E(n,b
+ h)}{n - J}},\quad g(n,b,h) = \frac{(n + J)(n -
b)}{n^3}\frac{1}{\frac{\tilde E(n,b)}{n - J}}\frac{1}{\frac{\tilde
E(n,b + h)}{n - J}}.$$
Set
$$\Delta _h^b\tilde E(n,b) = \mathcal {P}_h(n,b) =\mathcal {P},$$
then it follows
\begin{equation}\label{2-8}
\frac{d\mathcal {P}}{dn} \buildrel \Delta \over = f(n,b,h) -
g(n,b,h)\frac{\mathcal {P}}{n - J}.
\end{equation}
We also deduce from $\left. \tilde E(n,b) \right|_{n = J} = 0$ and
 $\left. \tilde E(n,b+h) \right|_{n = J} = 0$ that
\begin{equation}\label{2-9}
 \left.  \mathcal {P}  \right|_{n = J}  = \left.  \mathcal {P}_h (n,b)  \right|_{n = J} = 0.
\end{equation}
Now, we formally suppose that %
$$ \lim  \limits_{n \to J}
\frac{d\mathcal {P}}{dn} = l ,$$
 then by the Hospital's rule in form, we have
$$l = f(J,b,h) - g(J,b,h)l.$$
It is easy to see that $$ l = \frac{f(J,b,h)}{1 + g(J,b,h)} < 0,
$$%
and $$l \to  l_0 \buildrel \Delta \over = \frac{f(J,b)}{1 +
g(J,b)},\quad {\rm as} \quad h \to 0.$$

On the other hand, by \eqref{2-5} we obtain
\begin{equation}\label{2-10}
\begin{split} & f(n,b,h) \to f(J,b,h) < 0,\quad {\rm as} \quad n \to J, \\
 & g(n,b,h) \to g(J,b,h),\quad {\rm as} \quad n \to J,
\end{split}\end{equation}
and
\begin{equation*}\begin{split}
& f(n,b,h) \to f(n,b,0) < 0,\quad  {\rm as} \quad h \to 0, \\
 & g(n,b,h) \to g(n,b,0),\quad  {\rm as} \quad h \to 0, \\
& g(J,b,h) > C > 0,\quad g(n,b,h) > C > 0.
\end{split}\end{equation*}

Next we make singularity analysis on problem \eqref{2-8}-\eqref{2-8}
with respect to $n$ in a small neighborhood around the singularity
point $n=J$. With the help of the same methods in \cite{WMZZ21}, we
want to prove that $\frac{\mathcal {P}}{n-J}$ has not only an upper
bound but also a lower bound.

We claim that there exist two positive constants $\varepsilon_*$ and
$\delta_*$ independent of $b$ such that
$$ l -  \delta _*  < \frac{\mathcal {P}}{ n - J } < l +  \delta _*,
\quad \forall n \in [J-\varepsilon_*, J+\varepsilon_*].$$

\begin{figure}
\centering
\includegraphics[width=6cm]{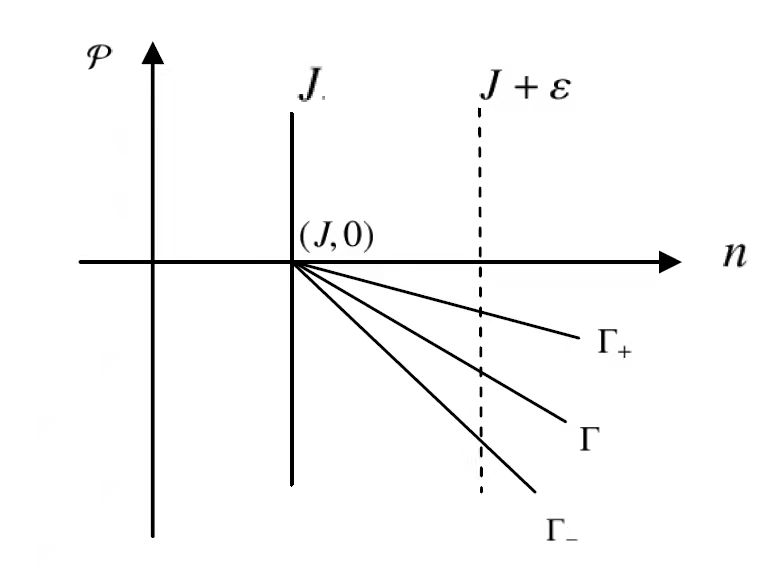}
\includegraphics[width=6cm]{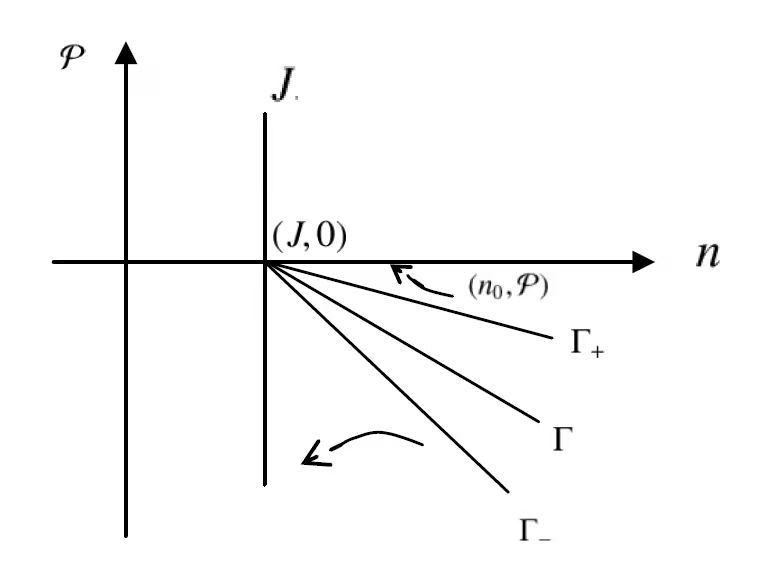}
\caption{Sketched graph  of $\mathcal {P}$ and $n$ near the singular point.}
\end{figure}



Indeed, we let $\Gamma_+, \Gamma$ and $\Gamma_-$ be three rays which
pass through the point $(J,0)$ with the slope values $l+\delta_*$,
$l$ and $l-\delta_*$, respectively. We use $\Sigma_{\varepsilon}$ to
stand for the triangle area bounded by rays $\Gamma_+, \Gamma_-$ and
the straight line $n=J+\varepsilon$ (see Figure 4). We would like to
show that there exists a positive constant $\varepsilon_+$ suitable
small such that $\mathcal {P}$ must be in $\Sigma_{\varepsilon_+}$
as $n\in [J, J+\varepsilon_+]$.

We use the method of proof by contradiction. For any $n_0\in
[J,J+\varepsilon]$, we assume that there always exists one point
$(n_0, \mathcal {P})$ which is lying above the triangle area
$\Sigma_{\varepsilon}$.

At $n=n_0$, it follows from \eqref{2-8} that
\begin{equation*}
\frac{d\mathcal {P}}{dn} = f(n,b,h) - g(n,b,h)\frac{\mathcal {P}}{n
- J}.
\end{equation*}
This together with  $ \frac{\mathcal {P}}{n - J} > l + \delta _*$
gives
\begin{equation*}\begin{split}
 \frac{{d\mathcal {P}}}{{dn}}({n_0}) =& f({n_0},b,h) - g({n_0},b,h)\frac{\mathcal {P}}{{{n_0} - J}} \\
  \le& f({n_0},b,h) - g({n_0},b,h)(l + \delta _*) \\
  \le& f(J,b,h) - g(J,b,h)l \\
  =& l < l + \delta _*.
 \end{split}\end{equation*}
This deduces that the trajectory pass through $(n_0, \mathcal
{P}(n_0))$ will intersect with the $n-$axis before $n=J$ (see Figure
5). However, this yields a contradiction with $\mathcal {P}(J)=0$.
Therefore, we obtain that there must exist constant
$\varepsilon_+>0$ such that $(n,\mathcal {P})\in
\Sigma_{\varepsilon_+}$ for any $n\in [J,J+\varepsilon_+]$.
Furthermore, %
$$\frac{\mathcal {P}}{n-J}\in [l-\delta_*, l+\delta_*], \quad \forall n
\in [J,J+\varepsilon_+].$$
Obviously, the above process holds independent of $b$ and provided
that $h>0$ is suitable small. Similarly, we can prove that  there
must exist constant $\varepsilon_->0$ suitable small such that
$(n,\mathcal {P})\in \Sigma_{\varepsilon_-}$ for any $n\in
[J-\varepsilon_-, J]$,
which implies, %
$\frac{\mathcal {P}}{n-J}\in [l-\delta_*, l+\delta_*], ~~ \forall n
\in [J-\varepsilon_-,J].$

By letting $$\varepsilon_*=\min\{\varepsilon_+, \varepsilon_-\},$$
we obtain
$$l-\delta_* \le \frac{\mathcal {P}}{n-J} \le l+\delta_*, \quad \forall n
\in [J-\varepsilon_*,J+\varepsilon_*].$$

On the other hand, since \eqref{2-8} has no singularity on the
domain $[n_*, J-\varepsilon_*]\cup [J+\varepsilon_*, M]$, it is
easy to get
$$\left|\frac{\mathcal {P}}{n-J}\right|<C, \quad \forall n\in
[n_*, J-\varepsilon_*]\cup [J+\varepsilon_*, M].$$
Hence, it follows that
$$\left|\frac{\mathcal
{P}}{n-J}\right|<C, \quad \forall n\in [n_*, M],$$
which implies that
$$\left|\frac{\Delta _h^b\tilde E(n,b)}{n-J}\right|<C, \quad \forall n\in [n_*, M],$$
namely,
$$\left| \frac{\tilde E(n,b + h)
- \tilde E(n,b)}{n-J}\right|<Ch, \quad \forall n\in [n_*, M],$$
which gives \eqref{2-6}.
The proof of Lemma \ref{L2-5} is completed. \hfill $\Box$

\medskip

Based on the primary works above, we begin to prove the stability of
the $C^1$-smooth transonic solution as follows.
\begin{lem}\label{L2-6}
There exists a constant $C>0$  such that
\begin{equation}
\label{(14)-8} \|  n _1(\cdot) -  n _2 (\cdot)\|_{ C^1 [0,L]}
  + \|  E_1(\cdot) -  E_2 (\cdot)\|_{ C^2 [0,L]} \le C\delta_0.
\end{equation}
\end{lem}
\noindent\textbf{Proof.} By plugging $(n, \tilde E, b)=(n_1, \tilde
E(n_1, b_1), b_1)$ and $(n, \tilde E, b)=(n_2, \tilde E(n_2, b_2),
b_2)$ into \eqref{2-3}, respectively. And taking difference of the
resulted equations, we have
\begin{equation*}
  \begin{split}
\left( n_1- n_2 \right)_x
 =& \frac{n_1^3}{n_1 + J}\frac{\tilde E(n_1,b_1)}{n_1 - J}
 - \frac{n_2^3}{n_2 + J}\frac{\tilde E(n_2,b_2)}{n_2 - J} \\
  =& \frac{n_1^3}{n_1 + J}\frac{\tilde E(n_1,b_1)}{n_1 - J}
  - \frac{n_2^3}{n_2 + J}\frac{\tilde E(n_1,b_1)}{n_1 - J}
  + \frac{n_2^3}{n_2 + J}\frac{\tilde E(n_1,b_1)}{n_1 - J}
   - \frac{n_2^3}{n_2 + J}\frac{\tilde E(n_2,b_2)}{n_2 - J} \\
  = &\frac{\tilde E(n_1,b_1)}{n_1 - J}
  \left( \frac{n_1^3}{n_1 + J}
  - \frac{n_2^3}{n_2 + J} \right)
  + T(n_1,n_2,b_1,b_2),
 \end{split}
\end{equation*}
where
\begin{equation*}\begin{split}
 T(n_1,n_2,b_1,b_2) =& \frac{n_2^3}{n_2 + J}
 \left( \frac{\tilde E(n_1,b_1)}{n_1 - J}
 - \frac{\tilde E(n_2,b_1)}{n_2 - J}
 + \frac{\tilde E(n_2,b_1)}{n_2 - J}
  - \frac{\tilde E(n_2,b_2)}{n_2 - J} \right).
 \end{split}\end{equation*}
By using Lemma \ref{L2-4} and the estimate \eqref{2-6} in Lemma
\ref{L2-5}, we have
\begin{equation*}\begin{split}
 T(n_1,n_2,b_1,b_2)
  =& \frac{n_2^3}{n_2 + J}
  \left( \frac{d\left(\frac{\tilde E(n_1,b_1)}{n_1 - J} \right)}{n}
  (n = \xi )(n_1 - n_2)
  + \frac{\tilde E(n_2,b_1)
   - \tilde E(n_2,b_2)}{(n_2 - J)(b_1 - b_2)}(b_1 - b_2) \right) \\
  \le& C\left| n_1 - n_2 \right| + C\left| b_1 - b_2
  \right|,
 \end{split}\end{equation*}
where $\xi$ is a number between $n_1$ and $n_2$.

Then, by Lemma \ref{L2-4} again, we have
\begin{equation*}
 \begin{split}
 \left| n_1 - n_2 \right|_x
   \le &\frac{\tilde E(n_1,b_1)}{n_1 - J}
  \left( \frac{n_1^3}{n_1 + J}
  - \frac{n_2^3}{n_2 + J} \right)+ C\left| n_1 - n_2 \right|
  + C\left| b_1 - b_2 \right|\\
  \le & C\left| n_1 - n_2 \right| + C\left| b_1 - b_2 \right|,
 \end{split}
    \end{equation*}
which implies
\begin{equation*}
\frac{d\left( n_1 - n_2 \right)^2}{dx} \le C\left( \left| n_1 - n_2
\right|^2 + \left|b_1 - b_2 \right|^2 \right),
\end{equation*}
where the Cauchy-Schwarz inequality was used.
Therefore, it follows that
\begin{equation*}
 \begin{split}
\left( n_1 - n_2 \right)^2 \le C\left( \left| n_{10} - n_{20}
\right|^2 + \left| b_1 - b_2 \right|^2 \right), \quad x\in [0,L],
 \end{split}
    \end{equation*}
i.e.
\begin{equation*}  \begin{split}
\left| n_1 - n_2 \right| \le C\left( \left| n_{10} - n_{20} \right|
+ \left| b_1 - b_2 \right|  \right)\quad x\in [0,L].
 \end{split}
    \end{equation*}

On the other hand, from the second equation of \eqref{2-0}, we
obtain
\begin{equation*} \begin{split}
  E_i (x)=E(n_i, b_i) =  E_{i0}  + \int_0^x
 \left(  n _i (y) - b_i  \right)dy, \quad i=1,2.\end{split} \end{equation*}
Then it follows  that \eqref{1-17}- \eqref{1-19}
are also true for $\alpha>0$. Namely, we prove
\[
  \|  n _1(\cdot) -  n _2 (\cdot)\|_{ C^1 [0,L]}
  + \|  E_1(\cdot) -  E_2 (\cdot)\|_{ C^2 [0,L]} \le C\delta_0.
 \]
The proof of Lemma \ref{L2-6} is completed. \hfill $\Box$

Finally, by combining Lemma \ref{L2-2}
 and
 Lemma \ref{L2-6}, we immediately prove the structural stability
 of  $C^1$-smooth transonic steady-states of \eqref{1.7} on $[0,L]$
 in Theorem \ref{T2}.
\section{Structural stability for steady transonic shock solutions}
In this section, we mainly prove Theorem \ref{T3} and establish the
structural stability for steady transonic shock solutions. Since Theorem \ref{T3-2} can be similarly obtained, we omit its proof.
\subsection{Preliminaries}
First, we prove the monotonic relation between the shock position
and the downstream density and a priori estimates for the steady
flows, which play a crucial role for the proof of Theorem \ref{T3}.
For any supersonic state $(n, E)$ satisfying $n < J$ , we can
connect it to a unique subsonic state $(\mathscr{S}(n), E)$ via a
transonic shock. Here $\mathscr{S}(n)$ is determined by the entropy
condition and the Rankine-Hugoniot condition
\begin{equation}\label{2.1}
\mathscr{S}(n) + \frac{{{J^2}}}{{\mathscr{S}(n)}} = n +
\frac{{{J^2}}}{n},\quad {\rm and} \quad \mathscr{S}(n) > J.
\end{equation}
By differentiating \eqref{2.1} with respect to $n$, we have
\begin{equation}\label{2.2}
\frac{{d\mathscr{S}(n)}}{{dn}} = \displaystyle\frac{{1 -
\displaystyle\frac{{{J^2}}}{{{n^2}}}}}{{1 -\displaystyle
\frac{{{J^2}}}{{{\mathscr{S}^2}(n)}}}}.
\end{equation}
This together with \eqref{1.6} give
\begin{equation}\label{2.3}
\displaystyle\frac{{d\mathscr{S}(n(x))}}{{dx}} = \frac{{nE - \alpha
J}}{{1 - \displaystyle\frac{{{J^2}}}{{{\mathscr{S}^2}(n)}}}}.
\end{equation}
\begin{lem}[Monotonic relation for the transonic shock solutions] \label{L2.1}

Let $(n^{(1)}, E^{(1)})$ and $(n^{(2)}, E^{(2)})$ be two transonic
shock solutions of \eqref{1.6}, and $(n^{(i)}, E^{(i)})$ $(i = 1,
2)$ are defined by
\begin{equation*} ( {{n^{(i)}},{E^{(i)}}} ) = \left\{
{\begin{array}{*{20}{c}}
 \vspace{3mm}  { ( {n_{\rm sup }^{(i)},E_{\rm sup }^{(i)}}  )(x),} \hfill
  & {{  as} \quad x\in (0, x_i),} \hfill  \\
   { ( {n_{\rm sub}^{(i)},E_{\rm sub}^{(i)}}  )(x),} \hfill &
   {{as} \quad x\in (x_i, L),} \hfill  \\
\end{array}} \right.\end{equation*}
where
\begin{equation*}
n_{\rm sup }^{(i)} < J < n_{\rm sub}^{(i)},\quad {  for} \quad i =
1,2.
\end{equation*}

They satisfy the same upstream boundary conditions
\begin{equation*}
{n^{(1)}}(0) = {n^{(2)}}(0) = {n_l}, \quad {E^{(1)}}(0) =
{E^{(2)}}(0) = {E_l}.\end{equation*}
Then, if $b < J, {x_1} < {x_2}$ and $E_{\rm sup }^{(2)}({x_1}) > 0$,
we have

$${n^{(1)}}(L) > {n^{(2)}}(L).$$
\end{lem}
\noindent\textbf{Proof.} For $x \in [0,{x_1})$, due to the fact that
both $  ( {n_{\sup }^{(1)},E_{\sup }^{(1)}} )$ and $ ( {n_{\sup
}^{(2)},E_{\sup }^{(2)}} ) $ satisfy  ordinary differential
equations \eqref{1.6} and the same initial data, we obtain
$$
 ( {n_{\sup }^{(1)},E_{\sup }^{(1)}} ) = ( {n_{\sup
}^{(2)},E_{\sup }^{(2)}} ), \quad {\rm as} \quad x \in [0,{x_1}).
$$

For $x\in [{x_1}, {x_2}]$, we define a function $\mathscr{E}$ as
follows
\begin{equation*}\left\{ {\begin{array}{*{20}{c}}
 \vspace{3mm}  {\displaystyle\frac{{d\mathscr{E}}}{{dx}} = \mathscr{S}
 ( {n_{\sup }^{(2)}} ) - b,
 \quad {\rm as} \quad x\in[x_1,x_2],} \hfill & {} \hfill  \\
   {\mathscr{E}({x_1}) = E_{\rm sub}^{(1)}({x_1}) = E_{\sup }^{(1)}({x_1}) = E_{\sup }^{(2)}({x_1}).} \hfill & {} \hfill  \\
\end{array}} \right.\end{equation*}

In view of $ n_{\rm sup }^{(2)} < J < \mathscr{S} ( {n_{\rm sup
}^{(2)}} )$, by using the comparison principles for  ordinary
differential equations \cite{Pa92}, we get $E_{\sup }^{(2)}(x) <
\mathscr{E}(x)$ as $x \in ({x_1},{x_2}].$ This, together with
\eqref{2.3}, gives
$$\frac{{d\mathscr{S} ( {n_{\sup }^{(2)}} )}}{{dx}} = \frac{{n_{\sup
}^{(2)}E_{\sup }^{(2)} - \alpha J}}{{1 -
\displaystyle\frac{{{J^2}}}{{{\mathscr{S}^2}(n_{\sup }^{(2)})}}}} <
\frac{{n_{\sup }^{(2)}\mathscr{E}(x) - \alpha J}}{{1 -
\displaystyle\frac{{{J^2}}}{{{\mathscr{S}^2}(n_{\sup }^{(2)})}}}}.$$
On the other hand, it follows from $ \mathscr{E}({x_1}) = E_{\sup
}^{(2)}({x_1})
> 0$ and $b < J < \mathscr{S} ( {n_{\rm sup }^{(2)}}
 ) $ that $\mathscr{E}(x) > 0$ as $x \in ({x_1},{x_2}]$, which
furthermore gives
\begin{equation*}\left\{ {\begin{array}{*{20}{c}}
  \vspace{3mm} {\displaystyle\frac{{d\mathscr{S} ( {n_{\rm sup }^{(2)}} )}}{{dx}}
  < \displaystyle\frac{{\mathscr{S} ( {n_{\rm sup }^{(2)}} )\mathscr{E}(x)
   - \alpha J }}{{1
    - \displaystyle\frac{{{J^2}}}{{{\mathscr{S}^2}(n_{\rm sup }^{(2)})}}}},} \hfill  \\
\vspace{3mm}   {\displaystyle\frac{{d\mathscr{E}}}{{dx}}
= \mathscr{S} ( {n_{\rm sup }^{(2)}}  ) - b,} \hfill  \\
   {\mathscr{S} ( {n_{\rm sup }^{(2)}}  )({x_1})
   = n_{\rm sub}^{(1)}({x_1}),\mathscr{E}({x_1}) = E_{\rm sub}^{(1)}({x_1}).} \hfill  \\
\end{array}} \right.\end{equation*}
Then, by using the comparison principles for  ordinary differential
equations again, we have
$$
\mathscr{S}\left( {n_{\rm sup }^{(2)}} \right)({x_2}) < n_{\rm
sub}^{(1)}({x_2}),\quad \mathscr{E}({x_2}) < E_{\rm
sub}^{(1)}({x_2}).$$
In view of $E_{\sup }^{(2)}({x_2}) < \mathscr{E}({x_2}),$ we obtain
$$
E_{\rm sub}^{(2)}({x_2}) = E_{\sup }^{(2)}({x_2}) <
\mathscr{E}({x_2}) < E_{\rm sub}^{(1)}({x_2}).$$

Recall that $ ( {n_{\rm sub}^{(1)},E_{\rm sub}^{(1)}}  )$ and $ (
{n_{\rm sub}^{(2)},E_{\rm sub}^{(2)}}  )$ solve the same  ordinary
differential equations on $[x_2, L]$, by the comparison principle
for  ordinary differential equations once more, we get
$$
n_{\rm sub}^{(1)}(L) > n_{\rm sub}^{(2)}(L),\quad {\rm and}\quad
E_{\rm sub}^{(1)}(L)
> E_{\rm sub}^{(2)}(L).
$$
 The proof of Lemma
\ref{L2.1} is completed. \hfill $\Box$

\vspace{3mm}

 Next, by using the multiplier method, we establish the a
priori estimates for supersonic and subsonic flows, which yield the
existence of supersonic, subsonic, and transonic shock solutions.

It follows from \eqref{1.6}  that $n$ satisfies
\begin{equation}\label{2.4}
\frac{d}{dx}\left( {\mathfrak{f}(n)n_x +
\alpha\displaystyle\frac{J}{{ n}}} \right) = n - b,
\end{equation}
in which $\mathfrak{f}(n) = \displaystyle \frac{n^2-J^2}{n^3}$.

Let $(n_0, E_0)$ be a supersonic or subsonic solution of \eqref{1.6}
with the doping profile $b_0$ and with initial data $(n_I , E_I )$,
i.e.,
\begin{equation}\label{2.5}
\left\{ {\begin{array}{*{20}{c}}
 \vspace{3mm} \displaystyle\frac{d}{dx}\left( {\mathfrak{f}(n)n_x +
\alpha\displaystyle\frac{J}{{ n}}} \right) = n - {b_0},  \hfill  \\
   {n(a) = {n_I},\quad n_x(a) = \displaystyle\frac{{{E_I}
   - \alpha\displaystyle\frac{J}{{{n_I}}}}}{{\mathfrak{f}({n_I})}}.} \hfill  \\
\end{array}} \right.
\end{equation}

In the following lemma, we give  the stability estimates for both
the supersonic and the subsonic solutions of \eqref{1.6}, which are
small perturbations of the solutions to the problem \eqref{2.5}.

\begin{lem}[] \label{L2.2}

For any interval $[a, l]\subseteq [0,x_0) \cup (x_0, L]$, suppose
$(n_0, E_0)$ to be a supersonic or subsonic solution to the problem
\eqref{2.5}. Then there is $\epsilon> 0$ such that if
\begin{equation}\label{2.6}
{\left\| {b(x) - {b_0}} \right\|_{{C^0}[a,l]}} + \left| {{{\tilde
n}_I}} \right| + \left| {{{\tilde E}_I}} \right| < \epsilon ,
\end{equation}
then,  for $x \in [a, l]$, there exists a unique supersonic or
subsonic solution $(n, E)(x)$ to the problem \eqref{1.6} with
initial conditions
\begin{equation}\label{2.7}
n(a) = {n_I} + {\tilde n_I}, \quad E(a) = {E_I} + {\tilde E_I}.
\end{equation}
Furthermore, $(n, E)$ satisfies
\begin{equation}\label{2.8}
{\left\| {n - {n_0}} \right\|_{{C^1}[a,l]}} < C{e^{\gamma
L}}\epsilon ,
\end{equation}
where constants $C > 0$ and $\gamma > 0$.
\end{lem}
\noindent\textbf{Proof.} The proof is given only for the case when
$n_0$ is supersonic on $[a, l]$ (the case when $n_0$ is subsonic is
quite similar). When $n_0$ is supersonic on $[a, l]$, there exist
constants $c_0 > 0$,  $c_1 > 0$ and $c_2 > 0$ such that
\begin{equation}\label{2.9}
\mathfrak{f}\left( {{n_0}} \right)(x) > {c_0}, \quad  {c_1}<
{n_0}(x) < J,\quad {\rm and} \quad \left| {\frac{d}{dx}n_0(x)}
\right| \le {c_2},\quad {\rm as} \quad x \in [a,l].
\end{equation}
First, we prove the results by assuming that
\begin{equation}\label{2.10}
\mathfrak{f}\left( n \right)(x) > \frac{{{c_0}}}{2}, \quad
\frac{1}{2}{c_1} < n(x) < J,\quad {\rm and} \quad \left|
\frac{d}{dx}n(x) \right| \le 2{c_2},\quad {\rm as} \quad x \in
[a,l].
\end{equation}

If we get the estimate \eqref{2.8}, then the lemma can be proved by
using the local existence theory of  ordinary differential equations
and the standard continuation argument.

Set $\tilde n = n - {n_0}$
 and $\tilde b = b - {b_0}$. Then from \eqref{2.4}-\eqref{2.5}, we obtain
\begin{equation}\label{2.11}
\frac{d}{dx}\left( {\mathfrak{f}(n)\tilde n_x +
{\mathfrak{F}_1}({n_0},\tilde n)\tilde n\frac{dn_0}{dx}   +
{\mathfrak{F}_2}({n_0},\tilde n)\tilde n} \right) - \tilde n =  -
\tilde b,
\end{equation}
where $${\mathfrak{F}_1}({n_0},\tilde n) = \int_0^L
\frac{d}{dn}\mathfrak{f} ({n_0} + \theta \tilde n)d\theta \quad {\rm
and} \quad {\mathfrak{F}_2}({n_0},\tilde n) = \alpha J\int_0^L
\displaystyle\frac{{ - 1}}{{{{\left( {{n_0} + \theta \tilde n}
\right)}^2}}}d\theta.$$

For constant $ \mu> 0$, we define a multiplier $\mathfrak {K}(x) :=
e^{-\mu(x-a)}$, and then multiply both sides of \eqref{2.11} by
$\mathfrak{K}(x)(\tilde n_x+\tilde n )$. By an integration by parts,
we have
\begin{equation*}\begin{split}
 & - \int_a^l   \tilde b \mathfrak{K}(x)(\tilde n_x + \tilde n)dx \\
  =& \int_a^l   \mathfrak{K}(x)
  \left( \left( {\frac{\mu }{2} - 1} \right)\mathfrak{f}(n)
  + \frac{1}{2} \frac{d}{dx}  \mathfrak{f}(n)    + {\mathfrak{F}_1}({n_0},
  \tilde n)\frac{dn_0}{dx}
   + {\mathfrak{F}_2}({n_0},\tilde n)\right) {(\tilde n_x)^2}dx \\
  &+ \int_a^l   \mathfrak{K}(x)\left( {\frac{{{\mu
  ^2}}}{2}\mathfrak{f}(n)
  - \frac{\mu }{2}\frac{d}{dx}  \mathfrak{f}(n)
   + \frac{d}{dx}\left( \mathfrak{F}_1 ({n_0},\tilde n)\frac{dn_0}{dx}\right) - 1
  +  \frac{d}{dx} \mathfrak{F}_2}({n_0},\tilde n)  \right){(\tilde n)^2}dx \\
  &+ \int_a^l   \mathfrak{K}(x)\left( {\frac{d}{dx}
  \left( \mathfrak{F}_1 ({n_0},\tilde n)\frac{dn_0}{dx}\right)
  +   \mathfrak{F}_1 ({n_0},\tilde n)\frac{dn_0}{dx}  - 1
   +   \mathfrak{F}_2 ({n_0},\tilde n)
  +   \frac{d}{dx} \mathfrak{F}_2 ({n_0},\tilde n)} \right)\tilde n\tilde n_xdx \\
  &+ \mathfrak{K}(l){\left[ {\mathfrak{f}(n)\frac{{{{(\tilde n_x)}^2}}}{2}
   + \mathfrak{f}(n)\tilde n\tilde n_x
  + \mu \mathfrak{f}(n)\frac{{{{(\tilde n)}^2}}}{2}} \right]_{x = l}}
   - {\left[ {\mathfrak{f}(n)\frac{{{{(\tilde n_x)}^2}}}{2}
   + \mathfrak{f}(n)\tilde n\tilde n_x
   + \mu \mathfrak{f}(n)\frac{{{{(\tilde n)}^2}}}{2}} \right]_{x = a}} .
 \end{split}\end{equation*}

It follows from \eqref{2.10} that we can choose $\mu$  large enough
such that
\begin{align}\label{2.12}
& \int_a^l  \mu {e^{ - \mu (x - a)}}\left( {\mu {{\tilde n}^2}
 + {{\left( {\tilde n_x} \right)}^2}} \right)dx + {e^{ - \mu (l - a)}}
 \left( {\mu {{\tilde n}^2} + {{\left( {\tilde n_x} \right)}^2}} \right)(l) \\
  \le & C\int_a^l   {e^{ - \mu (x - a)}}{{\tilde b}^2}dx
  + C\left( {\mu \tilde n_I^2 + \tilde E_I^2} \right). \notag
\end{align}

On the other hand, from \eqref{2.6}, we get
$$\int_a^l {} {e^{ - \mu (x - a)}}{\tilde b^2}dx + C\left(
{\mu \tilde n_I^2 + \tilde E_I^2} \right) \le C{\epsilon ^2}.$$
This, together with \eqref{2.12}, yields $$\int_a^l   \left(
{{{\tilde n}^2} + {{\left( {\tilde n_x} \right)}^2}} \right)dx \le
C{e^{\mu (l - a)}}{\epsilon ^2} \le C{e^{\mu }}{\epsilon ^2}.$$
By Sobolev embedding Theorem, it further gives
$$ \left\|  \tilde n
\right\|_{ C^0 [a,l]}  \le \left\|  \tilde n \right\|_{ H^1 [a,l]}
\le C{e^{\gamma L}}\epsilon,$$
 where $\gamma =\displaystyle
\frac{\mu }{2}.$  %
Then it follows from \eqref{2.11} that
$${\left\| {\tilde n_x} \right\|_{{C^0}[a,l]}} \le C{e^{\gamma L }}\epsilon.
$$
Hence, \eqref{2.8} follows by combining the above two inequalities.
The proof of the lemma \ref{L2.2} is completed.
 \hfill $\Box$

\subsection{Structural stability for transonic shock
solutions}
 In this subsection, we begin to prove
Theorem \ref{T3} and establish the structural stability for
transonic shock solutions to the boundary value problem \eqref{1.6}
and \eqref{1.4}.

\medskip

\noindent\textbf{Proof of Theorem \ref{T3}.} The proof is divided
into three stepsn.

\emph{Step 1.} For $b(x) = b_0(x)$ and  $i=1, 2$, we prove that
there exist transonic shock solutions $(n_i, E_i)(x)$ with the shock
location at $x_i$ such that $x_2
> x_1$, and
$n_2^r(L) < {n_r} < n_1^r(L)$.

By the conditions stated in Theorem \ref{T3} on
 the unperturbed transonic shock solution
$(n^{(0)}, E^{(0)})$ for the case when $b(x) = b_0(x) (x \in [0,
L])$, there is a constant $\delta > 0$ which satisfies $[x_0-\delta,
x_0 + \delta] \subset (0, L)$, such that the ordinary differential
equations
\begin{equation}\label{2.13}
\begin{split}
\frac{d}{dx}{\left( {n + \frac{{{J^2}}}{n}} \right)} = nE -\alpha
J,\quad {E_x} = n - {b_0}, \end{split}
\end{equation}
with the initial condition
\begin{equation}\label{2.14}
\begin{split}
 (  n,E  ) \left|   \right._{x = 0}  = (  n_l , E_l
 ), \end{split}
\end{equation}
have a unique smooth solution $(n^l, E^l)(x)$ on the interval $x \in
[0, x_0 + \delta]$ which satisfies $0 < n^l(x) < J$ for $x \in [0,
x_0 + \delta]$ and
\begin{equation}\label{2.15}
\begin{split}
 {E^l}(x) > 0,\quad {\rm for} \quad x\in [x_0 - \delta,x_0 + \delta] , \end{split}
\end{equation}
where $x_0$ is the shock location for
 $(n^{(0)}, E^{(0)})$ for the case when
$b(x) = b_0 (x \in [0, L])$. Furthermore, it follows from the
uniqueness for the initial value problems of ordinary differential
equations that
$$
 ( n^l , E^l )(x) =  (  n^{(0)} , E^{(0)} )(x),
\quad {\rm as} \quad x \in [0,{x_0}).$$

Set ${x_1} = {x_0} - \delta$ and ${x_2} = {x_0} + \delta .$
Then for $x \in [{x_i},L]$, let $\left( {n_i^r,E_i^r} \right)(x)(i =
1,2)$ be the solution of the ordinary differential equations
\eqref{2.13} with the initial conditions
$$
\ ( n_i^r,E_i^r  )|_{x=x_i} =  ( \mathscr{S} (
 n^l ( x_i )   ), E^l ( x_i )   ),\quad {\rm for} \quad i =
1, 2.
$$
We obtain from Lemma \ref{L2.2} that there is a unique smooth
subsonic solution $(n^r_i, E^r_i)(x)$ on the interval $x \in [x_i ,
L]$ satisfying $n^r_i (x) > n_s$ and $E^r_i (x) > E^l(x_i) > 0$ as
$x \in (x_i , L] (i = 1, 2)$. Furthermore, Lemma \ref{L2.1} together
with $x_2
> x_1$ yield
\begin{equation}\label{2.16}
\begin{split}
n_2^r(L) < {n_r} < n_1^r(L). \end{split}
\end{equation}

 \emph{Step 2.} When $b$ is a small perturbation of $b_0$, we prove
that there exist two transonic shock solutions $ (  \hat n ^{(1)},
\hat E ^{(1)} )(x)$ and $ ( \hat n ^{(2)}, \hat E^{(2)} )(x)$ with
the shock location at $x_1$ and $x_2$, respectively, such that $\hat
n_2^r(L) < {n_r} < \hat n_1^r( L )$.

For the case that $b$ is a small perturbation of $b_0$, we define
two transonic solutions based on $\left( {n_1^r,E_1^r} \right)$ and
$\left( {n_2^r,E_2^r} \right)$. Let
\begin{equation*} ( {{{\hat n}^{(i)}},{{\hat E}^{(i)}}}  )(x) = \left\{
{\begin{array}{*{20}{c}}
 \vspace{3mm}  { ( {\hat n_i^l,\hat E_i^l} )(x),} \hfill
 & {{\rm as} \quad x\in [0, x_i),} \hfill  \\
   { ( {\hat n_i^r,\hat E_i^r}  )(x),} \hfill
   & {{\rm as} \quad x\in (x_i, L],} \hfill  \\
\end{array}} \right.\end{equation*}
for $i=1,2,$ where $ ( {\hat n_i^l,\hat E_i^l}  )(x)$ is the
solution of the ordinary differential equations
\begin{equation}\label{2.17}
\begin{split}
\frac{d}{dx}{\left( {n + \frac{{{J^2}}}{n}} \right) } = nE -\alpha
J,\quad {E_x} = n - b(x), \end{split}
\end{equation}
on the region $[0, x_i ]$ with the initial data \eqref{2.14} and
$(\hat n^r_i, \hat E^r_i )$ is the solution of the ordinary
differential equations \eqref{2.17}
on $[x_i , L]$ with the initial data %
$$(\hat n^r_i, \hat E^r_i )|_{x=x_i}
=  ( {\mathscr{S} ( {\hat n_i^l({x_i} - )}  ),\hat E_i^l({x_i})}
 ).$$

By Lemma \ref{L2.2} and \eqref{1.10}, %
we obtain that $(\hat n_i^l, \hat E_i^l )$ and $(\hat n_i^r, \hat
E_i^r )$ are well-defined and satisfy
$$\mathop {\sup }\limits_{x \in [0,{x_i})} \left| { ( {\hat
n_i^l,\hat E_i^l}  ) - ( {n_i^l,{E^l}}  )} \right| \le C\varepsilon
,\quad \mathop {\sup }\limits_{x \in ({x_i},L]} \left| { ( {\hat
n_i^r,\hat E_i^r}  ) -  ( {n_i^r,E_i^r}  )} \right| \le C\varepsilon
, \quad {\rm for}\quad i = 1, 2.$$
 Furthermore, we have
$$
 \left| {\hat n_i^r(L) - n_i^r( L)} \right| \le
C\varepsilon , \quad {\rm for}\quad i = 1, 2.$$
This, together with \eqref{2.16}, yields that
$$\hat n_2^r(L) < {n_r} < \hat n_1^r( L),$$
provided that $\epsilon>0$ is small enough.

\emph{Step 3.} We prove that there exists a unique transonic shock
solution $(\tilde n, \tilde E  )$ with a single transonic shock
located at a point $\tilde{x}_0 \in (x_1, x_2)$.

The results established in steps 1-2 show that the boundary problem
\eqref{1.6} and \eqref{1.4} admits a unique transonic shock solution
$(\tilde n, \tilde E  )$ with a single transonic shock located at
some point $\tilde{x}_0 \in (x_1, x_2)$ by a monotonicity argument
as follows: for $x \in [x_1, x_2]$, we define a function
$\mathfrak{M}(x) = n(L)$ where $n$ is a transonic shock solution of
the system \eqref{1.6} satisfying \eqref{2.14} with shock located at
$x$. By Lemmas \ref{L2.1}-\ref{L2.2}, we obtain that
$\mathfrak{M}(x)$ is continuous strictly decreasing on $[x_1, x_2]$.
Furthermore, it follows from the stability estimate \eqref{2.8} in
Lemma \ref{L2.2} that $\tilde x_0 \in [x_0- C\epsilon , x_0 +
C\epsilon ]$. We have completed the proof of Theorem \ref{T3}.
 \hfill $\Box$

\section{Linear dynamic instability of transonic shock solutions}
In this section, we study the linear dynamic instability of
transonic shock solutions for the Euler-Poisson equations with
relaxation effect \eqref{1.1}.
\subsection{Formulation of the linearized problem}\label{s3.1}
Let $(\bar n, \bar u, \bar E)$ be a steady transonic shock solution
of the form \eqref{1.11} which satisfies \eqref{1.24}. Assume that
the initial values $(n_0, u_0, E_0)$ satisfy \eqref{1.20} and the
 compatibility conditions. From \cite{LY85}, we obtain that there
 is a piecewise smooth solution
containing a single shock $x = s(t)$ (with $s(0) = \tilde x_0)$
which satisfies the Rankine-Hugoniot conditions \eqref{1.23} and the
Lax geometric shock condition, of the Euler-Poisson equations with
relaxation effect on $[0, \bar T]$ for some $\bar T > 0$, which can
be written in the following form
\begin{equation}\label{3.1}
\begin{split}
(n,u,E)(t,x) = \left\{ {\begin{array}{*{20}{c}}
 \vspace{3mm}  {({n_ - },{u_ - },{E_ - }),} \hfill & {{\rm as}
 \quad x\in (0, s(t)),} \hfill  \\
   {({n_ + },{u_ + },{E_ + }),} \hfill & {{\rm as}
   \quad x\in (s(t), L).} \hfill  \\
\end{array}} \right.
\end{split}
\end{equation}
By noting that, when $t > T_0$ for some $T_0 > 0$, $(n_-, u_-, E_-)$
will depend only on the boundary values at $x = 0$. Furthermore,
when $\varepsilon$ is small, by the standard lifespan argument in
\cite{LY85}, we get $T_0 < \bar T$. Hence,
\begin{equation}\label{3.2}
\begin{split}
({n_ - },{u_ - },{E_ - }) = ({\bar n_ - },{\bar u_ - },{\bar E_ -
})\quad {\rm as} \quad t > {T_0}.
\end{split}
\end{equation}

In what follows,  we set $T_0 = 0$ for convenience. We would like to
extend the local-in-time solution to all $t > 0$. Note \eqref{3.2},
we need only to establish uniform estimates in the region $\{(t, x)|
t
> 0, x
> s(t)\}$. To this end, let us formulate an initial boundary value
problem in this region. Obviously, the Rankine-Hugoniot conditions
for \eqref{3.1} are
\begin{equation}\label{3.3}
\begin{split}
[nu] = [n]s'(t),\quad \left[ n +{n{u^2}} \right] = [nu]s'(t),
\end{split}
\end{equation}
where $[g] = g\left( {t ,s(t)^+} \right) - g\left( {t ,s(t)^-}
\right)$ is the jump of $g(t,x)$ at $x=s(t)$.

From \eqref{3.3}, we have $\left[ {n + n{u^2}} \right]\times[n] =
{\left[ {nu} \right]^2}. $
That is

$$
\left( {n_ + + \frac{{J_ + ^2}}{{{n_ + }}} - n_ - - \frac{{J_ -
^2}}{{{n_ - }}}} \right)\left( {{n_ + } - {n_ - }} \right)\left(
{t,s(t)} \right) = {\left( {{J_ + } - {J_ - }} \right)^2}\left(
t,s(t) \right),$$
where $J=nu$, $J_+=n_+u_+$ and  $J_-=n_-u_-$. By noting \eqref{3.2},
we get
\begin{equation}\label{3.4}
\begin{split}
J\left( {t,s(t)^- } \right) = \bar J.
\end{split}
\end{equation}
Therefore, we obtain
\begin{equation*}\begin{split}
&\left( {{n_ + }\left( {t,s(t)} \right) - {n_ - }\left( {s(t)} \right)} \right)\times\\
& \left( {n_ +\left( {t,s(t)} \right)
 \! + \!  \frac{{J_ + ^2\left(
{t,s(t)} \right)}}{{{n_ + }\left( {t,s(t)} \right)}}
 \!  - \!  {\bar n}_ +\left( {s(t)} \right)
 \!  -  \! \frac{{{{\bar J}^2}}}{{{{\bar n}_ + }}}\left( {s(t)} \right)}
 \!  + \!  {\bar n}_ +\left( {s(t)} \right)
  \! + \!  \frac{{{{\bar J}^2}}}{{{{\bar n}_ + }}}\left( {s(t)} \right)
   \! -  \! {\bar n}_
 -\left( {s(t)} \right)
 \!  - \!  \frac{\bar J^2}{{\bar n}_ - }\left( s(t) \right)\right) \\
  &= {\left( {{J_ + }\left( {t,s(t)} \right) - \bar J} \right)^2}.
 \end{split}\end{equation*}
By the Rankine-Hugoniot conditions \eqref{1.14} and Taylor
expansions, we have
\begin{equation*}\begin{split}
& \left( { {{n_ + }\left( {t,s(t)} \right)
 - {{\bar n}_ + }\left( {s(t)} \right)}
 - \bar u^2_ +\left( {s(t)} \right)\left( {{n_ + }
 \left( {t,s(t)} \right) - {{\bar n}_ + }\left( {s(t)} \right)} \right)} \right. \\
 &\hspace{5mm} + 2\bar u_ +\left( {s(t)}
 \right)\left( {{J_ + }\left( {t,s(t)} \right) - \bar J} \right)
 + {\partial _x}\left( {{\bar n}_ + + \frac{{{{\bar J}^2}}}{{{{\bar n}_ + }}}}
  \right)\left( {{x_0}} \right)\left( {s(t) - {x_0}} \right) \\
 &\hspace{5mm} \left. { - {\partial _x}\left( {{\bar n}_ -
 + \frac{{{{\bar J}^2}}}{{{{\bar n}_ - }}}} \right)\left( {{x_0}} \right)
 \left( {s(t) - {x_0}} \right) + {\mathcal{R}}} \right) \times \left( {{{\bar n}_ + }
 \left( {{x_0}} \right) - {{\bar n}_ - }\left( {{x_0}} \right) + {\mathcal{R}}} \right) \\
  &= {\left( {{J_ + }\left( {t,s(t)} \right) - \bar J} \right)^2},
\end{split}\end{equation*}
where
$${\mathcal{R}} = \mathcal {O}(1)\left( {{{\left( {{n_ + }\left(  t,s(t)  \right) -
{{\bar n}_ + }\left( {s(t)} \right)} \right)}^2} + {{\left( {{J_ +
}\left( {t,s(t)} \right) - \bar J} \right)}^2} + {{\left( {s(t) -
{x_0}} \right)}^2}} \right).$$

From now on, we usually use $\mathcal{R}$ to stand for those
quadratic terms with different $\mathcal {O}(1)$ coefficients. Then
it follows from the implicit function Theorem that
\begin{equation}\label{3.5}
\begin{split}
\left( {{J_ + } - \bar J} \right)\left( {t,s(t)} \right) = {\mathcal
{M}_1}\left( {\left( {{n_ + } - {{\bar n}_ + }} \right)\left(
t,s(t)  \right),s(t) - {x_0}} \right),
\end{split}
\end{equation}
where $\mathcal {M}_1$ satisfies
$${\mathcal {M}_1}\left( {0,0} \right) = 0,
\quad \displaystyle\frac{{\partial {\mathcal {M}_1}}}{{\partial
\left( {{n_ + } - {{\bar n}_ + }} \right)}} =  -
\displaystyle\frac{{1 - \bar u^2_ + }}{{ 2\bar u _ +}}({x_0}), \quad
\displaystyle\frac{{\partial {\mathcal {M}_1}}}{{\partial \left(
{s(t) - {x_0}} \right)}} =  - \displaystyle\frac{{\left( {{{\bar n}_
+ } - {{\bar n}_ - }} \right){{\bar E}_ + }}}{{ 2\bar u _
+}}({x_0}).$$
By plugging \eqref{3.5} into \eqref{3.3}$_1$, we get
\begin{equation}\label{3.6}
\begin{split}
s'(t) = {\mathcal {M}_2}\left( {\left( {{n_ + } - {{\bar n}_ + }}
\right)\left(  t,s(t)  \right),s(t) - {x_0}} \right),
\end{split}
\end{equation}
where $\mathcal {M}_2$ satisfies
 $${\mathcal {M}_2}\left( {0,0} \right) =
0, \quad \displaystyle\frac{{\partial {\mathcal {M}_2}}}{{\partial
\left( {{n_ + } - {{\bar n}_ + }} \right)}} =    \frac{ \bar u^2_ +
 -1 }{{2{{\bar u}_ + }\left( {{{\bar n}_ + } - {{\bar n}_ - }}
\right)}}({x_0}), \quad\displaystyle\frac{{\partial {\mathcal
{M}_2}}}{{\partial \left( {s(t) - {x_0}} \right)}} = -
\displaystyle\frac{{{{\bar E}_ + }}}{{2{{\bar u}_ + }}}({x_0}).$$

From \eqref{1.1}$_3$, we obtain %
$${E_ + }(t,x) = {E_l} + \int_0^{s(t)} {} \left( {{n_ - } - b}
\right)(z)dz + \int_{s(t)}^x   \left( {{n_ + } - b}
\right)(z)dz,\quad {\rm as}\quad  x\in (s(t), L] .$$
 It follows from  \eqref{1.1}$_1$ and the Rankine-Hugoniot conditions \eqref{3.3} that
$${\partial _t}{E_ + } = {n_l}{u_l} - {n_ + }{u_ + }(t,x)   = \bar J
- {J_ + }(t,x).$$
Let $\mathcal {V} = {E_ + }(t,x) - {\bar E_ + }(x),$ then
$$
 {\mathcal{V}_t} = \bar J - {J_ + }, \quad {\mathcal{V}_x} = {n_ + } - {\bar n_ + }.$$
Hence, from the momentum equation in the Euler-Poisson equations
with relaxation effect \eqref{1.1}, we have
$$
{\left( {{E_ + } - {{\bar E}_ + }} \right)_{tt}} + {\partial
_x}\left( {{\bar n}_ + + \frac{{{{\bar J}^2}}}{{{{\bar n}_ + }}} -
n_ + - \frac{{ J_ + ^2}}{{{n_ + }}}} \right) + {n_ + }{E_ + } -
{\bar n_ + }{\bar E_ + } -\displaystyle \alpha\left( {{J_ + } - \bar
J} \right)
 = 0.$$
 Then
\begin{equation}\label{3.7}
\begin{split}
{\mathcal{V}_{tt}} + {\partial _x}\left( {{\bar n}_ + + \frac{{\bar
J}^2} {\bar n}_ +  -  ({\bar n}_ + + {\mathcal{V}_x}) -
\frac{{{{\left( {\bar J - {\mathcal{V}_t}} \right)}^2}}}{{{{\bar n}_
+ } + {\mathcal{V}_x}}}} \right) + {\bar E_ + }{\mathcal{V}_x} +
{\bar n_ + }\mathcal{V}  + \mathcal{V}{\mathcal{V}_x}+
\alpha{\mathcal{V}_t} = 0.
\end{split}
\end{equation}

We set $\eta  = \left( {{\eta _0},{\eta _1}} \right) = (t,x)$,
${\partial _i} = \displaystyle\frac{\partial }{{{\partial
{\eta_i}}}}$ and ${\partial _{ij}} = \displaystyle\frac{{{\partial
^2}}}{{{\partial {\eta_i}}{\partial {\eta_j}}}}$ for $i,j = 0,1$.
Then \eqref{3.7} can be rewritten as follows
\begin{equation}\label{3.8}
\begin{split}
\sum\limits_{0 \le i,j \le 1}   {\bar \mu_{ij}}\left(
{x,{\mathcal{V}_t},{\mathcal{V}_x}} \right){\partial
_{ij}}\mathcal{V} + \sum\limits_{0 \le i \le 1} {} {\bar
\beta_i}\left( {x,{\mathcal{V}_t},{\mathcal{V}_x}} \right){\partial
_i}\mathcal{V} + \bar \zeta\left(
{x,{\mathcal{V}_t},{\mathcal{V}_x}} \right)\mathcal{V} = 0,
\end{split}
\end{equation}
where $\bar \mu_{i j}, \bar \beta_i$ and $\bar \zeta$ are smooth
functions of their variables, and satisfy
\begin{align}\label{3.9}
{\mathcal {L}_0}\mathcal{V} =& \sum\limits_{0 \le i,j \le 1}  {{\bar
\mu}_{ij}}\left( {x,0,0} \right){\partial _{ij}}\mathcal{V} +
\sum\limits_{0 \le i \le 1}
 {\bar \beta_i}\left( {x,0,0} \right){\partial _i}\mathcal{V}
 + \bar \zeta\left( {x,0,0} \right)\mathcal{V} \\
  =& {\mathcal{V}_{tt}} - {\partial _x}\left( {\left( {1
   - \bar u^2_ +} \right){\mathcal{V}_x}} \right)
    + {\partial _x}\left( {2 \bar u_ +{\mathcal{V}_t}} \right)
    + {{\bar E}_ + }{\mathcal{V}_x} + {{\bar n}_ + }\mathcal{V} +\alpha{\mathcal{V}_t}.\notag
\end{align}

Moreover, we rewrite the Rankine-Hugoniot conditions
\eqref{3.5}-\eqref{3.6} as
\begin{equation}\label{3.10}
\begin{split}
{\mathcal{V}_t} =  - {\mathcal {M}_1}\left( {{\mathcal{V}_x},s(t) -
{x_0}} \right),\end{split}
\end{equation}
and
\begin{equation}\label{3.11}
\begin{split}
 s'(t) = {\mathcal {M}_2}\left( {{\mathcal{V}_x},s(t) - {x_0}} \right),
\end{split}
\end{equation}
respectively. Furthermore, by a straightforward computation, we have
\begin{equation*}\begin{split}
 \mathcal{V}\left( {t,s(t)} \right) =& {E_ + }\left( {t,s(t)} \right) - {{\bar E}_ + }\left( {s(t)} \right) \\
  =& {{ E}_ - }\left( {t,s(t)} \right) - {{\bar E}_ + }\left( {s(t)} \right) \\
  =& {{\bar E}_ - }\left( {s(t)} \right) - {{\bar E}_ + }\left( {s(t)} \right) \\
  = &{{\bar E}_ - }\left( {s(t)} \right) - {{\bar E}_ - }\left( {{x_0}} \right) + {{\bar E}_ + }\left( {{x_0}} \right) - {{\bar E}_ + }\left( {s(t)} \right) \\
  =& \left( {{\partial _x}{{\bar E}_ - }\left( {{x_0}} \right)
  - {\partial _x}{{\bar E}_ + }\left( {{x_0}} \right)} \right)
  \left( {s(t) - {x_0}} \right) + {\mathcal{R}}.
\end{split}
\end{equation*}
This together with \eqref{1.1}$_3$ implies that
\begin{equation}\label{3.12}
\begin{split}
 s(t) - {x_0} = {\mathcal {M}_3}\left( {\mathcal{V}\left( {t,s(t)} \right)} \right),
\end{split}
\end{equation}
where
$$
{\mathcal {M}_3}\left( 0 \right) = 0, \quad {\rm and} \quad
\frac{{\partial {\mathcal {M}_3}}}{{\partial \mathcal{V}}} =
\frac{1}{{{{\bar n}_ - }({x_0}) - {{\bar n}_ + }({x_0})}}.$$

By combining \eqref{3.10} and \eqref{3.12}, we get
\begin{equation}\label{3.13}
\begin{split}
{\partial _t}\mathcal{V} = {\mathcal {M}_4}\left(
{{\mathcal{V}_x},\mathcal{V}} \right), \quad {\rm at}\quad x = s(t),
\end{split}
\end{equation}
where
$$
{\mathcal {M}_4}\left( {0,0} \right) = 0, \quad
\displaystyle\frac{{\partial {\mathcal {M}_4}}}{{\partial
{\mathcal{V}_x}}} =   \displaystyle\frac{{1  - \bar u_ +
^2}}{{2{{\bar u}_ + }}}({x_0}), \quad \displaystyle\frac{{\partial
{\mathcal {M}_4}}}{{\partial \mathcal{V}}} =
-\displaystyle\frac{{{{\bar E}_ + }}}{{2{{\bar u}_ + }}}({x_0}).$$
By noting that on the right boundary, $x = L$, $\mathcal{V}$
satisfies
\begin{equation}\label{3.14}
\begin{split}
{\partial _x}\mathcal{V} = 0,\quad {\rm at} \quad x = L.
\end{split}
\end{equation}
We would like to obtain uniform estimates for $\mathcal{V}$ and $s$
which satisfy \eqref{3.8} and \eqref{3.12}-\eqref{3.14}.

In order to to reformulate the problem to the fixed domain $[x_0,
L]$, let us introduce the transformation
$$
 t' = t, \quad  x' = \left( {L - {x_0}} \right)\frac{{x -
s(t)}}{{L - s(t)}} + {x_0}, \quad \sigma \left( { t'} \right) = s(t)
- {x_0},$$
and let
\begin{equation}\label{3.15}
\begin{split}
{p_1}\left( {x',\sigma } \right) = \frac{{L - x'}}{{L - {x_0} -
\sigma \left( { t'} \right)}}, \quad {p_2}\left( \sigma \right) =
\frac{{L - {x_0}}}{{L - {x_0} - \sigma \left( { t'} \right)}}.
\end{split}
\end{equation}
Then \eqref{3.7} turns into the following form
\begin{equation*}\begin{split}
& {\mathcal{V}_{ t' t'}} + {p_2}{\partial _{x'}}\left( {
\frac{{{{\bar J}^2}}}{{{{\bar n}_ + }}} -
   p_2 {\mathcal{V}_{ x'}} - \frac{{{{\left( {\bar J - {\mathcal{V}_{ t'}}
   + \sigma '\left( { t'} \right){p_1}{\mathcal{V}_{ x'}}} \right)}^2}}}
   {{{{\bar n}_ + }
    + {p_2}{\mathcal{V}_{x'}}}}} \right) - 2{p_1}\sigma '\left( { t'} \right)
 \mathcal{V}_{ x' t'} + {\left( {\sigma '\left( {t'} \right){p_1}}
\right)^2}\mathcal{V}_{ x'  x'}\\
&
 + {p_2}\mathcal{V}{\mathcal{V}_{ x'}} - 2\frac{{{{\left( {\sigma '\left( {
t'} \right)} \right)}^2}}p_1}{{1 - {x_0} - \sigma \left( { t'}
\right)}}  {\mathcal{V}_{ x'}} - \alpha\sigma '\left( { t'}
 \right){p_1}{\mathcal{V}_{ x'}} + {p_2}{{\bar E}_ + }{\mathcal{V}_{ x'}}
 + \alpha{\mathcal{V}_{t'}} + {{\bar n}_ + }\mathcal{V} \\
  = &{p_1} \sigma''\left( {t'}
 \right){\mathcal{V}_{x'}}.
\end{split}
\end{equation*}

By a direct calculation, the equation \eqref{3.12} turns into
\begin{equation}\label{3.16}
\begin{split}
\sigma \left( {t'} \right) = {\mathcal {M}_3}\left(
{\mathcal{V}\left( {t, x' = {x_0}} \right)} \right),
\end{split}
\end{equation}
and \eqref{3.11}  changes into
\begin{equation}\label{3.17}
\begin{split}
\frac{{d\sigma }}{{d  t'}} = {\mathcal {M}_2}\left( {{p_2}(\sigma
){\mathcal{V}_{ x'}},\sigma ( t')} \right).
\end{split}
\end{equation}
By using \eqref{3.17} to denote the quadratic terms for $\sigma$ in
terms of $\mathcal{V}$ , we get, at $ x' = x_0$,
\begin{equation*} \begin{split}
 \frac{{d\sigma }}{{d  t'}} = &{\mathcal {M}_2}\left( {0,0} \right)
 + \frac{{\partial {\mathcal {M}_2}}}{{\partial {\mathcal{V}_{ x}}}}
 {p_2(\sigma)\mathcal{V}_{ x'}}
 + \frac{{\partial {\mathcal {M}_2}}}{{\partial (s(t) - {x_0})}}\sigma ( t')
 + {\mathcal{R}} \\
  = &  \frac{{ \bar u_ + ^2-1 }}{{2\left( {{{\bar n}_ + }
   - {{\bar n}_ - }} \right){{\bar u}_ + }}}({x_0}){p_2(\sigma)\mathcal{V}_{ x'}}
    - \frac{{{{\bar E}_ + }}}{{2{{\bar u}_ + }}}({x_0})\sigma ( t') +
    {\mathcal{R}},
\end{split}
\end{equation*}
which further implies
\begin{equation}\label{3.18}
\begin{split}
\frac{{d\sigma }}{{d t'}} + \frac{{{{\bar E}_ + }}}{{2{{\bar u}_ +
}}}({x_0})\sigma  = {\mathcal {N}_2}\left( {{\mathcal{V}_{
x'}},\mathcal{V}} \right),
\end{split}
\end{equation}
where $\mathcal {N}_2$ satisfies
$$
\left| {{\mathcal {N}_2}\left( {{\mathcal{V}_{ x'}},\mathcal{V}}
\right) + \frac{{1 - \bar u_ + ^2}}{{2\left( {{{\bar n}_ + } -
{{\bar n}_ - }} \right){{\bar u}_ +
}}}({x_0}){p_2(\sigma)\mathcal{V}_{ x'}}} \right| =|\mathcal{R}|\le
C\left( {\mathcal{V}_{ x'}^2 + {\mathcal{V}^2}} \right).$$

Obviously, in view of \eqref{3.16} and \eqref{3.17}, both $\sigma$
and $\sigma'$ can be denoted in terms of $\mathcal{V}$ and its
derivatives at $x' = x_0$. Then, after handling \eqref{3.13} with
\eqref{3.16} and \eqref{3.17}, we obtain
\begin{equation}\label{3.19}
\begin{split}
{\mathcal{V}_{ t'}} = {\mathcal {N}_1}\left( {{\mathcal{V}_{
x'}},\mathcal{V}} \right),\quad {\rm at} \quad x' = {x_0}.
\end{split}
\end{equation}
Equivalently, with the help of the implicit function Theorem once
more, we get
\begin{equation}\label{3.20}
\begin{split}
{\mathcal{V}_{x'}} = {\mathcal {N}_3}\left( {{\mathcal{V}_{
t'}},\mathcal{V}} \right),\quad {\rm at} \quad x' = {x_0},
\end{split}
\end{equation}
where 
\begin{equation*} \begin{split}
\mathcal {N}_3\left( {{\mathcal{V}_{t'}},\mathcal{V}} \right) =&
{\mathcal {N}_3}\left( {0,0} \right)
  + \frac{{\partial {\mathcal {N}_3}}}{{\partial {\mathcal{V}_{ t'}}}}{\mathcal{V}_{ t'}}
   + \frac{{\partial {\mathcal {N}_3}}}{{\partial \mathcal{V}}}\mathcal{V} + {\mathcal{R}} \\
  =& \frac{{2{{\bar u}_ + }}}{{1 - \bar u_ + ^2}}
  ({x_0}){\mathcal{V}_{t'}} + \frac{{{{\bar E}_ + }}}{{1 - \bar u_ + ^2}}({x_0})
  \mathcal{V} + {\mathcal{R}}.
\end{split}
\end{equation*}
Or, equivalently,
$$\left| {{\mathcal {N}_3}\left( {{\mathcal{V}_{t'}},\mathcal{V}} \right)
- \frac{{2{{\bar u}_ + }}}{{1 - \bar u_ +
^2}}({x_0}){\mathcal{V}_{t'}} - \frac{{{{\bar E}_ + }}}{{1 - \bar u_
+ ^2}}({x_0})\mathcal{V}} \right| =|\mathcal{R}|\le C\left(
{\mathcal{V}_{t'}^2 + {\mathcal{V}^2}} \right).$$

Next, we still use $x$ and $t$ to denote $ x'$ and $t'$,
respectively, for convenience. The problem becomes
\begin{equation}\label{3.21}
\begin{cases}
\mathcal {L}(x,\mathcal{V},\sigma )\mathcal{V} = \sigma
''(t){p_1}{\partial _x}\mathcal{V}, \quad
(t,x) \in \R_+ \times [{x_0},L], \\
 {\partial _x}\mathcal{V} = {{\xi}_1}\left(
{{\mathcal{V}_t},\mathcal{V}} \right){\mathcal{V}_t} +
{\omega_1}\left( {{\mathcal{V}_t},\mathcal{V}} \right)\mathcal{V},
\quad {\rm at} \quad x = {x_0}, \\
 {\partial _x}\mathcal{V} = 0,\quad {\rm at} \quad x = L, \\
 \sigma (t) = {\mathcal {M}_3}\left( {\mathcal{V}(t,{x_0})}
\right),
\end{cases}
\end{equation}
where, by using $\eta_0$ and $\eta_1$ to stand for $t$ and $x$,
respectively,

\begin{equation*}\begin{split}
\mathcal {L}(x,\mathcal{V},\sigma )\mathcal{U}  = & \sum\limits_{i,j
= 0}^1 {} {\mu_{ij}} \left( {x,\mathcal{V},\nabla \mathcal{V},\sigma
,\sigma '} \right){\partial _{ij}}\mathcal{U}
 + \sum\limits_{i = 0}^1 {} {\beta_i}\left( {x,\mathcal{V},\nabla \mathcal{V},\sigma ,\sigma '}
  \right)
 {\partial _i}\mathcal{U}\\
  &+ \zeta\left( {x,\mathcal{V},\nabla \mathcal{V},\sigma ,\sigma '} \right)\mathcal{U},
 \end{split}\end{equation*}
with
$${{\xi}_1}\left( {{\mathcal{V}_t},\mathcal{V}} \right) = \int_0^1 {} \frac{{\partial
{\mathcal {N}_3}}}{{\partial {\mathcal{V}_t}}}\left( {\theta
{\mathcal{V}_t},\theta \mathcal{V}} \right)d\theta , \quad
{\omega_1}\left( {{\mathcal{V}_t},\mathcal{V}} \right) = \int_0^1 {}
\frac{{\partial {\mathcal {N}_3}}}{{\partial \mathcal{V}}}\left(
{\theta {\mathcal{V}_t},\theta \mathcal{V}} \right)d\theta .$$
Moreover, we have $\mathcal {L}(x, 0, 0)\mathcal{U} = \mathcal
{L}_0\mathcal{U}$, and
\begin{equation}\label{3.22}
\begin{cases}
{\mu_{00}}\left( {x,\mathcal{V},\nabla \mathcal{V},\sigma ,\sigma '}
\right) = 1,\quad {\mu_{11}}(x,0,0,0,0) = \bar u^2_ + - 1,
\\
 {\mu_{01}}(x,0,0,0,0) = {\mu_{10}}(x,0,0,0,0) = {{\bar u}_ + },\\
 {\beta_0}(x,0,0,0,0) = \alpha + 2{\partial _x} \left( \bar u _ +
\right),\quad {\beta_1}(x,0,0,0,0)
 = {\partial _x}\left( {\bar u^2_ + } \right) + {{\bar E}_ + }, \\
 \zeta(x,0,0,0,0) = {{\bar n}_ + },\quad {{\xi}_1}(0,0) =
\frac{{2{{\bar u}_ + }}} {{1 - \bar u_ + ^2}}({x_0}),\quad
{\omega_1}(0,0) = \frac{{{{\bar E}_ + }}}{{1 - \bar u_ +
^2}}({x_0}).
\end{cases}
\end{equation}
Then the linearized problem is
\begin{equation}\label{3.23}
\begin{cases}
 \mathcal {L}(x,0,0)\mathcal{V} = 0,  \quad
(t,x) \in \R_+ \times [{x_0},L], \\
 {\partial _x}\mathcal{V} = \frac{{2{{\bar u}_ + }}}{{1 - \bar u_ +
^2}}({x_0}){\mathcal{V}_t} + \frac{{{{\bar E}_ + }}}{{1 - \bar u_ +
^2}}
({x_0})\mathcal{V}\quad {\rm at} \quad x = {x_0}, \\
 {\partial _x}\mathcal{V} = 0,\quad {\rm at} \quad x = L \\
 \mathcal{V}(0,x) = {\mathfrak{h}_1}(x), \quad {\mathcal{V}_t}(0,x)
= {\mathfrak{h}_2}(x), \quad x\in ({x_0}, L).
\end{cases}
\end{equation}
\subsection{Linear dynamic instability}\label{s3.2}
Let  $x_0\in [0, L)$ be the shock location for the steady transonic
shock solution, we investigate the linear dynamic instability for
the steady transonic shock solutions when $\bar E (x_0) < -\delta$,
where
 $\delta>0$ is a constant. We rewrite the linearized problem
\eqref{3.23} as
\begin{equation}\label{3.24}
\begin{cases}
 {\mathcal{V}_{tt}}  -  {\partial _x} \left(
{\left( {1   -   \bar u_ + ^2} \right){\mathcal{V}_x}} \right)  +
2{\partial _x}  \left( {{{\bar u}_ + }{\mathcal{V}_t}}
 \right)    +   {{\bar n}_ + }\mathcal{V}   +    {{\bar E}_ + }{\mathcal{V}_x}
 +   \alpha{\mathcal{V}_t}  = 0,
\quad
(t,x)\in\R_+\times({x_0},L), \\
 {\mathcal{V}_t} = \frac{{1 - \bar u_ + ^2}} {{2{{\bar u}_ +
}}}({x_0}){\mathcal{V}_x} - \frac{{{{\bar E}_ + }}}{{2{{\bar u}_ +
}}}
({x_0})\mathcal{V}, \quad {\rm at} \quad x = {x_0}, \\
 {\partial _x}\mathcal{V} = 0,\quad {\rm at} \quad x =
L.
\end{cases}
\end{equation}
It follows from \eqref{1.24} that
\begin{equation}\label{3.25}
\begin{split} {\bar E_ + }\left( {{x_0}} \right) <  -
\delta.\end{split}
\end{equation}

In order to prove the linear instability, we look for solutions to
the problem \eqref{3.24} of the form $\mathcal{V}(t,x) = e^{\nu t}
\mathcal{U}(x)$. A direct computation gives
\begin{equation}\label{3.26}
\begin{cases}
 \left( {1  \!   -   \!   \bar u_ + ^2}
\right) {\mathcal{U}_{xx}}
     -  \!   \left( {{\partial _x} \left(  { \bar u_ + ^2} \right)
    +    2\nu {{\bar u}_ + }     + \!  {{\bar E}_ + }}
    \!\right){\mathcal{U}_x}
   \!   - \! \left( \!{{\nu ^2}  \!
   + \!  2\nu {\partial _x}{{\bar u}_ + }
  \!\! - \!  \alpha\nu  \!  + \!  {{\bar n}_ + }} \right)
  \mathcal{U}  \!\! =  \!
  0,~~
(t,x) \!    \in   \R_+ \! \times  \!  ( x_0,\!  L), \\
 {\mathcal{U}_x} = \frac{{2{{\bar u}_ + }}}{{1
 - \bar u_ + ^2}}({x_0})
 \left( {\frac{{{{\bar E}_ + }}}{{2{{\bar u}_ + }}}
 ({x_0}) + \nu } \right)\mathcal{U},\quad {\rm at} \quad x = {x_0}, \\
 {\mathcal{U}_x} = 0, \quad {\rm at} \quad x = L.
 \end{cases}
\end{equation}
For a fixed parameter $\mathcal{U}(x_0) = \gamma > 0$, let us
consider
\begin{equation}\label{3.27}
\begin{cases}
 \left( {1 \! - \!  \bar u_ + ^2}
\right){\mathcal{U}_{xx}}
  \! -  \! \left( {{\partial _x}\left( { \bar u_ + ^2} \right)
  \! + \!  2\nu {{\bar u}_ + }  \! + \!  {{\bar E}_ + }} \right)
  {\mathcal{U}_x}
   \! -  \! \left( {{\nu ^2}  \! +  \! 2\nu {\partial _x}{{\bar u}_ + }
  \! - \!  \alpha\nu  \!  + \!  {{\bar n}_ + }} \right)\mathcal{U}
   \! =  \! 0, ~~ {\rm as}~~ x \in (x_0, +\infty), \\
  \mathcal{U}({x_0}) = \gamma  > 0,\quad {\mathcal{U}_x}({x_0})
  = \frac{{2{{\bar u}_ + }}}{{1 - \bar u_ + ^2}}({x_0})
 \left( {\frac{{{{\bar E}_ + }}}{{2{{\bar u}_ + }}}
 ({x_0}) + \nu } \right)\mathcal{U}.
\end{cases}
\end{equation}

By noting that  $0<\bar u_+ <1$ and \eqref{3.25}, it follows that if
$\nu = 0$, then $\mathcal{U}_x (x_0) < 0$. Hence, there exists
$\varsigma_1
> x_0$ such that $\mathcal{U}_x (x) < 0$ as $x\in [x_0, \varsigma_1]$.
On the other hand, if $\nu= -{\frac{{{{\bar E}_ + }}}{{{{\bar u}_ +
}}}({x_0})} $, then $\mathcal{U}_x (x_0) > 0$. We obtain that there
is $\varsigma_2
> x_0$ such that $\mathcal{U}_x (x)
> 0$ for $x\in [x_0, \varsigma_2]$.

Let
 $\varsigma = \min\{\varsigma_1, \varsigma_2\}.$
By the continuous dependence of the  ordinary differential equations
with respect to the initial data and the parameters, there exists a
number $ \nu\in (0, -{\frac{{{{\bar E}_ + }}}{{{{\bar u}_ +
}}}({x_0})})$ such that the problem \eqref{3.27} has a solution
$\mathcal{U}=\mathcal{U}(x)$ satisfying $\mathcal{U}_x (\varsigma) =
0$ which is a solution of \eqref{3.26} on $[x_0, \varsigma]$. This
shows that the linearized problem \eqref{3.24} or \eqref{3.23} can
have exponentially growing solutions. The proof of Theorem \ref{T4}
has been finished.
  \hfill $\Box$

\vspace{3mm}

\noindent {\sc Acknowledgments:}  This work was done when Y. Feng
visited McGill University supported by China Scholarship Council
(CSC) for the senior visiting scholar program (202006545001). He
would like to express his sincere thanks for the hospitality of
McGill University and CSC. The research of M. Mei was supported
by NSERC grant RGPIN 354724-2016.  The research of G. Zhang was supported by
 NSF of China (No. 11871012)


\end{document}